\documentclass[a4paper]{article}
\usepackage{amsfonts}
\usepackage{amssymb}
\usepackage{amsmath}
\usepackage{latexsym}
\usepackage{graphicx}
\newtheorem{thm}{Theorem}[section]
\newtheorem{cor}[thm]{Corollary}
\newtheorem{lem}[thm]{Lemma}
\newtheorem{ex}[thm]{Example}
\newtheorem{rem}[thm]{Remark}

\newcommand{\proof}{\medskip \noindent {\bf Proof. \ \ }}
\newcommand{\qed}{\null\hfill $\Box\;\;$ \medskip}

\begin{document}

\parbox{1mm}

\begin{center}
{\bf {\sc \Large Wavelets from Laguerre Polynomials and
Toeplitz-type Operators}}
\end{center}

\vskip 12pt

\begin{center}
{\bf Ondrej HUTN\'IK}
\end{center}


\hspace{5mm}\parbox[t]{10cm}{\fontsize{9pt}{0.1in}\selectfont\noindent{\bf
Abstract.} We study Toeplitz-type operators with respect to
specific wavelets whose Fourier transforms are related to Laguerre
polynomials. On the one hand, this choice of wavelets underlines
the fact that these operators acting on wavelet subspaces share
many properties with the classical Toeplitz operators acting on
the Bergman spaces. On the other hand, it enables to study
poly-Bergman spaces and Toeplitz operators acting on them from a
different perspective. Restricting to symbols depending only on
vertical variable in the upper half-plane of the complex plane
these operators are unitarily equivalent to a multiplication
operator with a certain function. Since this function is
responsible for many interesting features of these Toeplitz-type
operators and their algebras, we investigate its behavior in more
detail. As a by-product we obtain an interesting observation about
the asymptotic behavior of true polyanalytic Bergman spaces.
Isomorphisms between the Calder\'on-Toeplitz operator algebras and
functional algebras are described and their consequences are
discussed. } \vskip 24pt

\section{Introduction}

The structure of the space of wavelet transforms inside
$L_2(G,\mathrm{d}\nu)$ (the space of all square-integrable
functions on the affine group $G$ with respect to the left
invariant Haar measure $\mathrm{d}\nu$) was described in our
paper~\cite{hutnik3}. The key tool in this description is the
(Bargmann-type) transform giving an isometrical isomorphism of the
space $L_2(G,\mathrm{d}\nu)$ under which the space of wavelet
transforms is mapped onto tensor product of $L_2$-spaces where one
of them is the rank-one space generated by a suitable function.
This isomorphism is constructed according to the general
Vasilevski scheme of decomposition of Hilbert spaces presented
in~\cite{vasilevski} which was systematically applied in his
book~\cite{vasilevskibook}, and provides an interesting tool to
study Toeplitz-type operators based on the Calder\'on reproducing
formula and acting on wavelet subspaces. This idea was developed
in general setting in paper~\cite{hutnik2} where some results on
Toeplitz-type operators related to wavelets were also given. Then
this technique was used in~\cite{hutnik} for fundamental case
where the wavelet subspaces arise from functions whose Fourier
transforms are related to Laguerre polynomials. Note that in this
case the rank-one space is generated by Laguerre functions
$\ell_k(x)=\mathrm{e}^{-x/2} L_k(x)$ forming an orthonormal basis
in $L_2(\mathbb{R}_+)$, where $L_k(x)$ stands for the Laguerre
polynomial of degree $k$ and type 0 at $x$.

This paper is a further continuation of the above mentioned
research. Here we are interested in some properties of these
Toeplitz-type operators acting on wavelet subspaces in detail. In
accordance with paper~\cite{hutnik} we consider specific wavelets
which enable us to obtain an interesting generalization of the
classical case of Toeplitz operators acting on the Bergman space.
This choice of wavelets also underlines the intriguing patterns
which the corresponding wavelet subspaces and poly-Bergman spaces
share. In fact, the main source of motivation for our study comes
from the Vasilevski book~\cite{vasilevskibook} summarizing results
of author and his collaborators on classical Toeplitz operators
and their algebras on the Bergman space. The presented methods and
techniques in that book are very useful and inspiring for our
purposes. Thus, on the one hand, many obtained results of this
paper as well as paper~\cite{hutnikII} may be simply viewed as
"wavelet analogs" of results known for the classical Toeplitz
operators. On the other hand, the described generalization of
Toeplitz operators goes by another direction as in the case of
Toeplitz operators acting on weighted Bergman spaces, and gives
rise to many new and interesting questions. The organization of
this paper is as follows.

\paragraph{Section~\ref{section2}: Representation of wavelet subspaces} We
introduce basic notions and basic tools of our study including a
parameterized family of wavelets
$\{\psi^{(k)}\}_{k\in\mathbb{Z}_+}$ whose Fourier transform is
related to Laguerre polynomials $L_k(x)$. Then we introduce the
isometrical isomorphism
$$U: L_2(G,\mathrm{d}\nu) \to L_2(\mathbb{R})\otimes L_2(\mathbb{R}_+)$$
describing the structure of wavelet subspaces $A^{(k)}$ (related
to wavelet $\psi^{(k)}$ and Hardy-space functions) inside
$L_2(G,\mathrm{d}\nu)$. This result also allows us to describe the
direct and natural connection between wavelet subspaces and Hardy
spaces. Further, we construct operators $R_k$ and $R_k^*$
providing decomposition of orthogonal projection $P^{(k)}$ from
$L_2(G,\mathrm{d}\nu)$ onto wavelet subspace $A^{(k)}$ and of
identity operator $I$ on $L_2(\mathbb{R}_+)$. Note that the
operator $R_k$ is an exact analog of the Bargmann transform
mapping the Fock space $F_2(\mathbb{C}^n)$ onto
$L_2(\mathbb{R}^n)$. For more details on Bargmann transform in
connection with time-frequency analysis and Toeplitz-type
operators in this context see~\cite[Section 2]{huthut2}.

\paragraph{Section~\ref{section3}: Calder\'on-Toeplitz operators} The
representation of wavelet subspaces $A^{(k)}$ is then used to
study Toeplitz-type operators based on the Calder\'on reproducing
formula. Given a function (symbol) $a=a(\zeta)$, $\zeta\in G$, the
so-called Calder\'on-Toeplitz operator $T_a^{(k)}$ acting on
wavelet subspace $A^{(k)}$ is defined as follows $$T_a^{(k)} f =
P^{(k)} (a f), \quad f\in A^{(k)}.$$ The operators $T_a^{(k)}$
naturally arise in connection with localization in time-scale
analysis context and they are a filtered version (with $a$ being a
filter) of considered signal $f$ analyzed using wavelets
$\psi^{(k)}$. The main idea here is to use the constructed unitary
operators $R_k$ and $R_k^*$ to reduce the Calder\'on-Toeplitz
operator $T_a^{(k)}$ acting on wavelet subspace $A^{(k)}$ to its
unitary equivalent image $R_k T_a^{(k)} R_k^*$ acting on
$L_2(\mathbb{R}_+)$. Indeed, the main result of Theorem~\ref{CTO1}
which gives an easy and direct access to the properties of
Calder\'on-Toeplitz operators, then reads as follows:

\vskip 5pt \textsf{The Calder\'on-Toeplitz operator $T^{(k)}_{a}$
with a symbol $a=a(v)$, $v \in \mathbb{R}_+$, acting on wavelet
subspace $A^{(k)}$ is unitarily equivalent to the multiplication
operator $\mathfrak{A}_a^{(k)}=\gamma_{a,k}I$ acting on
$L_2(\mathbb{R}_+)$, where the function $\gamma_{a,k}:
\mathbb{R}_+\to \mathbb{C}$ has the form}
\begin{equation}\label{gamma0}
\gamma_{a,k}(\xi)=\int_{\mathbb{R}_{+}}
a\left(\frac{v}{2\xi}\right)\ell_k^2(v)\,\mathrm{d}v, \quad
\xi\in\mathbb{R}_+.
\end{equation}
As it can be seen, the function $\gamma_{a,k}$ is obtained by
integrating a dilation of a symbol $a=a(v)$ of a
Calder\'on-Toeplitz operator $T_a^{(k)}$ against a Laguerre
function of order $k$. This result extends the result of
Vasilevski for the classical Toeplitz operators acting on the
Bergman space (i.e., the case $k=0$ in our notation) in very
interesting way which differs from the case of Toeplitz operators
acting on weighted Bergman spaces studied in paper~\cite{GKV} and
then summarized in Vasilevski book~\cite{vasilevskibook}.
Moreover, in Section~\ref{section3} we present a number of results
following immediately from this equivalency including the
spectral-type representation of Calder\'on-Toeplitz operators
whose symbols depend only on imaginary coordinate in the upper
half-plane, as well as formulas for the Wick symbols and the star
product in terms of our function $\gamma_{a,k}$. In this section
we also briefly mention the form of operators $R_k T_a^{(k)}
R_k^*$ for a symbol $a$ depending only on horizontal variable in
the upper half-plane, and for a symbol $a$ in a special and
important product form. This presents a way how certain
pseudo-differential operators naturally appear in this context of
Calder\'on-Toeplitz operators, and opens the door to many
important and interesting questions.

\paragraph{Section~\ref{section4}: Some properties of function $\gamma_{a,k}$} Since
the function $\gamma_{a,k}$ given by~(\ref{gamma0}) is responsible
for many interesting features and behavior of the corresponding
Calder\'on-Toeplitz operator $T_a^{(k)}$ acting on $A^{(k)}$, in
this section we investigate its basic properties. We study the
limit behavior of higher order derivatives of $\gamma_{a,k}$ for
bounded as well as integrable symbols $a$ which provides
information about certain interesting properties of
$\gamma_{a,k}$. Then we give sufficient conditions for
$\gamma_{a,k}$ to be continuous on the whole $[0,+\infty]$: a
question which is closely related to the behavior of a symbol
$a=a(v)$, $v\in\mathbb{R}_+$, at a neighborhood of points $0$ and
$+\infty$. More precisely, we prove the following result:

\vskip 5pt \textsf{If $a=a(v)$ is a bounded symbol on
$\mathbb{R}_+$ such that the limits $\lim\limits_{v\to 0}
a(v)=a_0$ and $\lim\limits_{v\to+\infty} a(v)=a_\infty$ exist,
then for each $k\in\mathbb{Z}_+$ the function $\gamma_{a,k}$
belongs to the algebra $C[0,+\infty]$ of continuous functions on
$[0,+\infty]$. Moreover, $\gamma_{a,k}(+\infty)=a_0$ and
$\gamma_{a,k}(0)=a_\infty$ for each $k\in\mathbb{Z}_+$.}

\vskip 5pt \noindent The result states that the limit at infinity
and at zero of the function $\gamma_{a,k}$ \textit{is independent
of parameter $k$}. In fact, it depends only on a limit of the
corresponding symbol $a$, but not on the particularly chosen
Laguerre functions. This is surprising because the wavelet
transforms with Laguerre functions of order $k$ live, up to a
multiplier isomorphism, in the so-called true polyanalytic Bergman
space of order $k$, which is rather different from the classical
Bergman space of analytic functions (corresponding to $k=0$ in our
notation). Thus, this result contains the remarkable observation
that, asymptotically, all the true polyanalytic Bergman spaces
have the same behavior, when such behavior is observed via the
association with the Calder\'on-Toeplitz operator described in
Theorem~\ref{CTO1}. This result has some important consequences in
quantum physics, signal analysis and in the asymptotic theory of
random matrices, which are not yet completely understood.

\paragraph{Section~\ref{section5}: Isomorphism between the Calder\'on-Toeplitz
operator algebra and functional algebra} In connection with the
above mentioned results some operator algebras and a functional
dependence of Calder\'on-Toeplitz operators are described in this
section. In particular, for any $a(v)=\chi_{[0,\lambda]}(v)$ with
$\lambda\in \mathbb{R}_+$ (here $\chi_{[0,\lambda]}(v)$ is the
characteristic function of the interval $[0,\lambda]$) we have
$$\gamma_{\chi_{[0,\lambda]},k}(\xi) = 1-N_{2k}(2\lambda\xi)\,\mathrm{e}^{-2\lambda \xi}, \quad
\xi\in\mathbb{R}_+,$$ where $N_{2k}$ is a certain polynomial of
degree $2k$ (see the interesting result of Lemma~\ref{lemmaell} in
Appendix). Thus,
$$T_{\chi_{[0,\lambda]}}^{(k)} = R_k^*
\gamma_{\chi_{[0,\lambda]},k} R_k,$$ and each $C^*$-algebra
generated by the Calder\'on-Toeplitz operators
$T_{\chi_{[0,\lambda]}}^{(k)}$ is isometric and isomorphic to the
algebra $C[0,+\infty]$. Moreover, if we consider a symbol $a=a(v)$
from the algebra $L^{\{0,+\infty\}}_\infty(\mathbb{R}_+)$ of
bounded functions on $\mathbb{R}_+$ having limits at the endpoints
$0$ and $+\infty$ such that the corresponding function
$\gamma_{a,k}$ separates the points of $[0,+\infty]$, then we
construct the function
$$\nabla_{a,\lambda}^{(k)}(x) = -\frac{1}{\lambda}\ln(1-x)\int_{\mathbb{R}_+}
a(v)(1-x)^{v/\lambda}
L_k^2\left(-\frac{v}{\lambda}\ln(1-x)\right)\,\mathrm{d}v$$ on
$[0,1]$ such that each Calder\'on-Toeplitz operator $T_a^{(k)}$
acting on $A^{(k)}$ with a symbol $a=a(v)\in
L^{\{0,+\infty\}}_\infty(\mathbb{R}_+)$ is the function of
Toeplitz operator acting on the Bergman space
$A^{(0)}=\mathcal{A}^2(\Pi)$ with symbol $\chi_{[0,\lambda]}(v)$,
i.e.,
$$T_{a}^{(k)} = \nabla_{a,\lambda}^{(k)}\left(T_{\chi_{[0,\lambda]}}^{(0)}\right).$$
Interpretation and applicability of these results from the
viewpoint of localization in the time-frequency analysis are
discussed in the end of Section~\ref{section5}.

\section{Representation of wavelet subspaces}\label{section2}

In this section we summarize basic methods, constructions and
results from our previous works, mainly from~\cite{hutnik}. We use
the obvious notations: $\mathbb{R}$ ($\mathbb{C}$, $\mathbb{N}$)
are the sets of all real (complex, natural) numbers,
$\overline{\mathbb{R}}=\mathbb{R}\cup\{-\infty,+\infty\}$ is the
two-point compactification of $\mathbb{R}$, and $\mathbb{R}_{+}$
($\mathbb{R}_{-}$) are the positive (negative) half-lines with
$\chi_+$ ($\chi_-$) their characteristic functions, respectively.

\paragraph{The affine group and admissible wavelets on the real line} Let
$G=\{\zeta=(u,v);\,\, u\in \mathbb{R}, v>0\}$ be the locally
compact the so-called \textit{``$ax+b$''-group} with the group law
$$(u,v)\diamond(u',v')=(vu'+u,vv')$$ and the left-invariant Haar measure
$\mathrm{d}\nu(\zeta)=v^{-2}\,\mathrm{d}u\,\mathrm{d}v$. We
identify the group $G$ with the upper half-plane
$\Pi=\{\zeta=u+\mathrm{i} v; \,u\in\mathbb{R}, v>0\}$ in the
complex plane $\mathbb{C}$ (with $\mathrm{i}^2=-1$) equipped with
the hyperbolic measure $\mathrm{d}\nu$. Then
$L_2(G,\mathrm{d}\nu)$ is the space of all square-integrable
functions on $G$ with respect to the measure $\mathrm{d}\nu$. In
what follows $\langle\cdot, \cdot\rangle$ always means the inner
product on $L_2(\mathbb{R})$, whereas $\langle\cdot,
\cdot\rangle_G$ denotes the inner product on $L_2(G,
\mathrm{d}\nu)$. For their corresponding norms we use the notation
$\|\cdot\|$ and $\|\cdot\|_G$, respectively.

A function $\psi \in L_2(\mathbb{R})$ is called an
\textit{admissible wavelet} if it satisfies the admissibility
condition
\begin{equation}\label{admissibility}
\int_{\mathbb{R}_{+}}
|\hat{\psi}(x\xi)|^{2}\,\frac{\mathrm{d}\xi}{\xi}=1
\end{equation}for almost every $x\in\mathbb{R}$, where $\hat{\psi}$ stands
for the unitary Fourier transform $\mathcal{F}: L_2(\mathbb{R})\to
L_2(\mathbb{R})$ given by
$$\mathcal{F}\{g\}(\xi)=\hat{g}(\xi)=\int_{\mathbb{R}} g(x)\,\mathrm{e}^{-2\pi
\mathrm{i}x\xi}\,\mathrm{d}x.$$ Let us mention the interesting
constructions generalizing admissible wavelets using
group-theoretical methods which are given in recent papers of
Kisil, see~\cite{kisil1} and~\cite{kisil2}.

\paragraph{Wavelets related to Laguerre functions} In many problems in mathematical physics,
signal analysis, parameter identification, analytical and computer
studies, etc., the \textit{Laguerre functions}
$$\ell_n^{(\alpha)}(y) = \left[\frac{n!}{\Gamma(n+\alpha+1)}\right]^{1/2}
y^{\alpha/2}\,\mathrm{e}^{-y/2}L_n^{(\alpha)}(y), \quad y\in
\mathbb{R}_+,$$ naturally appear. Here, $\Gamma$ is the Euler
Gamma function, and $L_n^{(\alpha)}$ is the \textit{Laguerre
polynomial} of degree $n\in\mathbb{Z}_+=\mathbb{N}\cup\{0\}$ and
type $\alpha$ (may be considered as an arbitrary complex number)
given by
\begin{equation}\label{defLaguerre}
L_n^{(\alpha)}(y)=\frac{y^{-\alpha}\,\mathrm{e}^y}{n!}
\frac{d^n}{dy^n}\left(\mathrm{e}^{-y}y^{n+\alpha}\right)=\sum_{i=0}^{n}(-1)^{i}{n+\alpha\choose
n-i} \frac{y^i}{i!}, \quad y\in \mathbb{R}_+,
\end{equation}cf.~\cite[formula 8.970.1]{GR}. Recall that the system of functions
$\{\ell_n^{(\alpha)}(y)\}_{n\in\mathbb{Z}_+}$  for $\alpha>-1$
forms an orthonormal basis in the space $L_2(\mathbb{R}_+,
\mathrm{d}y)$, i.e.,
$$\int_{\mathbb{R}_+} \ell^{(\alpha)}_m(y)\,\ell^{(\alpha)}_n(y)\,\mathrm{d}y =
\delta_{mn}, \quad m,n\in\mathbb{Z}_+.$$ For $\alpha=0$ we will
simply write $L_n(y)$ and $\ell_n(y)$. In what follows we proceed
as in~\cite{hutnik}, and for $k\in\mathbb{Z}_+$ we consider the
functions $\psi^{(k)}$ and $\bar{\psi}^{(k)}$ on $\mathbb{R}$
defined on the Fourier transform side as follows
$$\hat{\psi}^{(k)}(\xi)=\chi_+(\xi)\sqrt{2\xi}\,\ell_k(2\xi),\quad
\textrm{and}
\quad\hat{\bar{\psi}}^{(k)}(\xi)=\hat{\psi}^{(k)}(-\xi),$$
respectively. Note that Ingrid Daubechies used $\varphi$ with
$\hat{\varphi}(\xi) = \sqrt{2\xi}\, \hat{\psi}^{(0)}(\xi)$ as the
wavelet function in her famous book~\cite{daubechies}. It is
obvious that for each $k\in\mathbb{Z}_+$ the functions
$\psi^{(k)}$ and $\bar{\psi}^{(k)}$ satisfy the
admissibility condition~(\ref{admissibility}). 
Then according to the Calder\'on reproducing formula,
cf.~\cite{calderon},
$$f(u)=\int_{\mathbb{R}_+}\left(\psi^{(k)}_v*\psi^{(k)}_v* f\right)(u)\,\frac{dv}{v^2},
\quad
g(u)=\int_{\mathbb{R}_+}\left(\bar{\psi}^{(k)}_v*\bar{\psi}^{(k)}_v*
g\right)(u)\,\frac{dv}{v^2},$$ for all $f\in H_2^+(\mathbb{R})$
and $g\in H_2^-(\mathbb{R})$, where $H_2^+(\mathbb{R})$, resp.
$H_2^-(\mathbb{R})$, are the Hardy spaces, i.e.,
\begin{align*}
H_2^+(\mathbb{R}) & = \Bigl\{h\in L_2(\mathbb{R}); \,\,\textrm{supp}\, \hat{h}\subseteq [0,+\infty)\Bigr\}; \\
H_2^-(\mathbb{R}) & = \Bigl\{h\in L_2(\mathbb{R});
\,\,\textrm{supp}\, \hat{h}\subseteq (-\infty,0]\Bigr\},
\end{align*}respectively. Here $*$ denotes the usual convolution on
$L_2(\mathbb{R})$, and $\psi_v(u)=v^{-1/2}\psi(u/v)$, $(u,v)\in
G$, is a dilation of $\psi$ on $L_2(\mathbb{R})$. It is well-known
that $H_2^+(\mathbb{R})$ and $H_2^-(\mathbb{R})$ are the only
proper invariant subspaces under the (quasi-regular)
representation $\rho$ of $G$ on $L_2(\mathbb{R})$ given by
$$(\rho_\zeta f)(x) =
\frac{1}{\sqrt{v}}\,f\left(\frac{x-u}{v}\right), \quad
\zeta=(u,v)\in G.$$


\begin{figure}
\begin{center}
\includegraphics[scale=0.9]{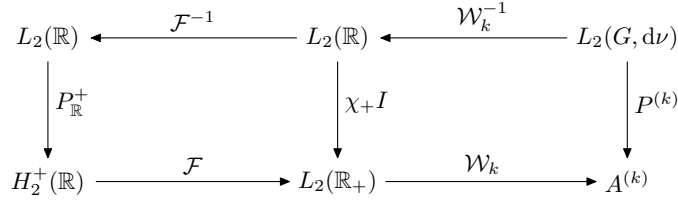}
\end{center}\caption{Relationship among the introduced spaces and operators}\label{fig_operatorsI}
\end{figure}

\paragraph{Continuous wavelet transform and wavelet subspaces} For
each $k\in\mathbb{Z}_+$ define the subspaces $A^{(k)}$ and
$\bar{A}^{(k)}$ of $L_2(G,\mathrm{d}\nu)$ as follows
\begin{align*}
A^{(k)} & := \left\{ (\mathcal{W}_k f)(u,v)=\left(f*\psi_v^{(k)}\right)(u); \,f\in H_2^+(\mathbb{R})\right\}; \\
\bar{A}^{(k)} & := \left\{(\mathcal{W}_{\bar{k}} f)(u,v) =
\left(f*\bar{\psi}_v^{(k)}\right)(u); \,f\in
H_2^-(\mathbb{R})\right\},
\end{align*}respectively. Note that $\mathcal{W}_k f$, resp. $\mathcal{W}_{\bar{k}} f$, are known as the
\textit{continuous wavelet transforms} of functions $f\in
H_2^+(\mathbb{R})$, resp. $f\in H_2^-(\mathbb{R})$, with respect
to wavelets $\psi^{(k)}$, resp. $\bar{\psi}^{(k)}$. Moreover,
$\mathcal{W}_k$ and $\mathcal{W}_{\bar{k}}$ are isometries from
$L_2(\mathbb{R})$ to $L_2(G, \mathrm{d}\nu)$ for each
$k\in\mathbb{Z}_+$. Consequently, $A^{(k)}$, resp.
$\bar{A}^{(k)}$, will be called the \textit{spaces of Calder\'on}
(or wavelet) \textit{transforms}. We also use the term wavelet
subspaces (of $L_2(G, \mathrm{d}\nu)$).

The relationship among the introduced spaces $A^{(k)}$ of wavelet
transforms of $H_2^+(\mathbb{R})$-functions, and the unitary
operators of continuous wavelet transform $\mathcal{W}_k$ and the
Fourier transform $\mathcal{F}$ is schematically described on
Figure~\ref{fig_operatorsI}. Note that $\mathcal{F}^{-1}:
L_2(\mathbb{R})\to L_2(\mathbb{R})$ and $\mathcal{W}_k^{-1}:
A^{(k)}\to L_2(\mathbb{R})$ are the inverse Fourier transform and
the inverse continuous wavelet transform, respectively. Observe
that Figure~\ref{fig_operatorsI} also includes the following
well-known classical result.

\begin{lem}\label{lem1}
The Fourier transform $\mathcal{F}$ gives an isometrical
isomorphism of the space $L_2(\mathbb{R})$ onto itself under which
\begin{itemize}
\item[(i)] the Hardy space $H_2^+(\mathbb{R})$, resp.
$H_2^-(\mathbb{R})$, is mapped onto $L_2(\mathbb{R}_+)$, resp.
$L_2(\mathbb{R}_-)$; i.e., $$\mathcal{F}: H_2^+(\mathbb{R}) \to
L_2(\mathbb{R}_+),\quad \textrm{resp}.\quad \mathcal{F}:
H_2^-(\mathbb{R}) \to
L_2(\mathbb{R}_-);$$  
\item[(ii)] the Szeg\"{o} projection $P_\mathbb{R}^+:
L_2(\mathbb{R})\to H_2^+(\mathbb{R})$, resp. $P_\mathbb{R}^-:
L_2(\mathbb{R})\to H_2^-(\mathbb{R})$, is unitarily equivalent to
the following one
$$\mathcal{F}P_\mathbb{R}^+\mathcal{F}^{-1} = \chi_+ I, \quad
\textrm{resp.}\quad \mathcal{F}P_\mathbb{R}^-\mathcal{F}^{-1} =
\chi_- I.$$
\end{itemize}
\end{lem}

For each $k\in\mathbb{Z}_+$ the spaces $A^{(k)}$ and
$\bar{A}^{(k)}$ are the reproducing kernel Hilbert spaces.
Explicit formulas for their reproducing kernels
$$K_\zeta^{(k)}(\eta)=\left\langle \rho_\eta\psi^{(k)},
\rho_\zeta\psi^{(k)}\right\rangle, \quad \textrm{resp.} \quad
\bar{K}_\zeta^{(k)}(\eta)=\overline{K_\zeta^{(k)}(\eta)},$$ and
orthogonal projections $P^{(k)}: L_2(G, \mathrm{d}\nu)\to
A^{(k)}$, resp. $\bar{P}^{(k)}: L_2(G, \mathrm{d}\nu)\to
\bar{A}^{(k)}$, are given in~\cite{hutnik}. Note also, that the
functions $\psi^{(k)}$ are not normalized in $L_2(\mathbb{R})$,
since
\begin{equation}\label{lambda_k}
\|\rho_\zeta\psi^{(k)}\|^2 = \|\psi^{(k)}\|^2 =
\|\hat{\psi}^{(k)}\|^2 = \frac{1}{2} \int_{\mathbb{R}_+} x\,
\ell_k^2(x)\,\mathrm{d}x = \frac{2k+1}{2} := \kappa_k,
\end{equation} where the formula
$$\int_{\mathbb{R}_+} x\,
\ell_k^2(x)\,\mathrm{d}x = 2k+1,\quad k\in\mathbb{Z}_+,
$$ has been used (as a special case of
formula~(\ref{Laguerreformula}) from Appendix). Therefore, the
constant $\kappa_k$ appears in formulas of Wick calculus, see
Section~\ref{section3}.

\paragraph{Structural results for wavelet subspaces} In what follows we introduce some important
operators used in our study, see~\cite{hutnik} for more details.
Interpret the space $L_2(G, \mathrm{d}\nu)$ as tensor product in
the form
$$L_2(G, \mathrm{d}\nu(\zeta)) = L_{2}(\mathbb{R}, \mathrm{d}u)
\otimes L_{2}(\mathbb{R}_{+}, v^{-2}\mathrm{d}v)$$ with $\zeta =
(u,v) \in G$, and consider the unitary operator $$U_{1}=
(\mathcal{F}\otimes I): L_{2}(\mathbb{R}, \mathrm{d}u) \otimes
L_{2}(\mathbb{R}_{+}, v^{-2}\mathrm{d}v) \to L_{2}(\mathbb{R},
\mathrm{d}\omega) \otimes L_{2}(\mathbb{R}_{+},
v^{-2}\mathrm{d}v).$$ For the purpose to "linearize" the
hyperbolic measure $\mathrm{d}\nu$ onto the usual Lebes\-gue plane
measure we introduce the unitary operator $$U_{2}:
L_{2}(\mathbb{R},\,\mathrm{d}\omega)\otimes
L_{2}(\mathbb{R}_{+},\,v^{-2}\mathrm{d}v) \to
L_{2}(\mathbb{R},\,\mathrm{d}x)\otimes
L_{2}(\mathbb{R}_{+},\,\mathrm{d}y)$$ given by
$$U_{2}: F(\omega,v) \longmapsto
\frac{\sqrt{2|x|}}{y} F\left(x,\frac{y}{2|x|}\right).$$  Then the
inverse operator
$$U_{2}^{-1}=U_{2}^{*}: L_{2}(\mathbb{R},\,\mathrm{d}x) \otimes
L_{2}(\mathbb{R}_{+},\,\mathrm{d}y) \to
L_{2}(\mathbb{R},\,\mathrm{d}\omega) \otimes
L_{2}(\mathbb{R}_{+},\,v^{-2}\mathrm{d}v)$$ is given by the rule
$$U_{2}^{-1}: F(x,y) \longmapsto \sqrt{2|\omega|}\,v F(\omega,2|\omega|v).$$
Using the classical result of Lemma~\ref{lem1} and the operator
$U=U_2 U_1$ we get the following theorem describing the structure
of wavelet subspaces $A^{(k)}$ and $\bar{A}^{(k)}$ inside
$L_2(G,\mathrm{d}\nu)$, see~\cite[Theorem~2.1]{hutnik} for its
proof.

\begin{thm}\label{thm1}
The unitary operator $U=U_{2}U_{1}$ gives an isometrical
isomorphism of the space $L_{2}(G,\mathrm{d}\nu)$ onto
$L_{2}(\mathbb{R},\,\mathrm{d}x)\otimes
L_{2}(\mathbb{R}_{+},\,\mathrm{d}y)$ under which
\begin{itemize}
\item[(i)] the space $A^{(k)}$ is mapped onto $L_2(\mathbb{R}_+)
\otimes L_{k}$, where $L_{k}$ is the rank-one space generated by
Laguerre function $\ell_{k}(y)=\mathrm{e}^{-y/2}L_k(y)$;
\item[(ii)] the projection $P^{(k)}: L_2(G,\mathrm{d}\nu)\to
A^{(k)}$ is unitarily equivalent to the following one
$$UP^{(k)}U^{-1} = \chi_+ I\otimes P_0^{(k)},$$ where $P_0^{(k)}$
given by $$\left(P_0^{(k)}H\right)(y)= \ell_k(y)
\int_{\mathbb{R}_+} H(t)\ell_k(t)\,\mathrm{d}t$$ is the
one-dimensional projection of $L_2(\mathbb{R}_+, \mathrm{d}y)$
onto $L_k$.
\end{itemize}
\end{thm}

The theorem may be stated analogously for the space
$\bar{A}^{(k)}$. Moreover, we may say more about the connection
between the wavelet subspaces and Hardy spaces. Indeed, as a
direct consequence of Lemma~\ref{lem1} and Theorem~\ref{thm1} we
have

\begin{thm}\label{thm2}
The unitary operator $V=(\mathcal{F}^{-1}\otimes
I)U_{2}(\mathcal{F}\otimes I)$ gives an isometrical isomorphism of
the space $L_{2}(G,\mathrm{d}\nu)$ onto
$L_{2}(\mathbb{R},\,\mathrm{d}x)\otimes
L_{2}(\mathbb{R}_{+},\,\mathrm{d}y)$ under which
\begin{itemize}
\item[(i)] the spaces $A^{(k)}$ and $H_2^+(\mathbb{R})$ are
connected by
the formula $$V\left(A^{(k)}\right) = H_2^+(\mathbb{R})\otimes L_k;$$  
\item[(ii)] the projections $P^{(k)}$ and $P_\mathbb{R}^+$ are
connected by the formula
$$VP^{(k)}V^{-1} = P^+_\mathbb{R}\otimes P_0^{(k)}.$$
\end{itemize}
\end{thm}


\begin{figure}
\begin{center}
\includegraphics[scale=0.9]{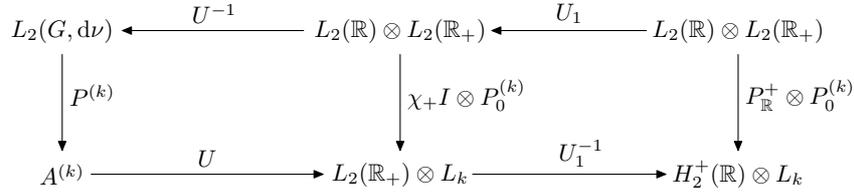}
\end{center}\caption{Visualizing the results of Theorem~\ref{thm1}
and Theorem~\ref{thm2}}\label{fig_operatorsII}
\end{figure}

The analogous result holds for $\bar{A}^{(k)}$,
$H_2^-(\mathbb{R})$, and $\bar{P}^{(k)}$, $P_\mathbb{R}^-$,
respectively. The diagram on Figure~\ref{fig_operatorsII}
schematically describes all the relations among the constructed
operators and spaces appearing in the above two theorems. Note
that the constructed operators $U_1$ and $U_2$ may serve also for
other purposes, e.g., may be useful in the study of certain
operator algebras, but here we will not continue in this
direction.

\begin{rem}\rm
The connection between spaces of wavelet transforms (with respect
to the specific Bergman wavelet and functions from the Hardy
space) and Bergman spaces is well-known,
see~\cite[Theorem~3.1]{hutnik2} and references given therein. In
fact, this result was a source of motivation for our research
in~\cite{hutnik2} for general wavelets. On the other hand, the
connection between the Hardy spaces and poly-Bergman spaces was
described in~\cite[Theorem~4.5]{vasilevski2}. From this point of
view the above results reveal that poly-Bergman spaces and wavelet
subspaces share intriguing patterns that may prove usable. For the
deeper study of this connection see the paper~\cite{abreu} with
some interesting applications to wavelet (super)frames. Moreover,
this suggested technique was recently successfully used to obtain
a complete characterization of all lattice sampling and
interpolating sequences in the Bargmann-Fock space of polyanalytic
functions, cf.~\cite{abreu2}, having a great potential in various
applications.
\end{rem}


\begin{figure}
\begin{center}
\includegraphics[scale=0.9]{operators3.1}
\end{center}\caption{Decomposition of orthogonal projection $P^{(k)}$
and identity $I$ on $L_2(\mathbb{R}_+)$}\label{fig_operatorsIII}
\end{figure}

\paragraph{Construction of Bargmann-type transforms} For the purpose
to construct an exact analog of the Bargmann transform (mapping
the Fock space $F_2(\mathbb{C}^n)$ onto $L_2(\mathbb{R}^n)$),
first let us introduce the isometric imbedding
$$Q_k: L_{2}(\mathbb{R}_+) \to
L_{2}(\mathbb{R})\otimes L_{2}(\mathbb{R}_{+})$$ given by
$$\left(Q_k f\right)(x,y) = \chi_+(x) f(x)\ell_{k}(y).$$ Here the function $f$
is extended to an element of $L_2(\mathbb{R})$ by setting
$f(x)\equiv 0$ for $x<0$. 
Its adjoint operator
$$Q_k^{*}: L_{2}(\mathbb{R})\otimes
L_{2}(\mathbb{R}_{+}) \to L_{2}(\mathbb{R}_+)$$ is given by
$$\left(Q_k^{*} F\right)(x) = \chi_{+}(x) \int_{\mathbb{R}_{+}}
F(x,t)\ell_{k}(t)\,\mathrm{d}t.$$ Then the operator $R_{k}=Q_k^{*}
U$ maps the space $L_{2}(G,\mathrm{d}\nu)$ onto
$L_{2}(\mathbb{R}_+)$, and the restriction
$$R_{k}|_{A^{(k)}}: A^{(k)} \to L_2(\mathbb{R}_+)$$ is
an isometrical isomorphism. The adjoint operator
$$R_{k}^{*}=U^{*} Q_k: L_{2}(\mathbb{R}_+) \to
A^{(k)} \subset L_{2}(G,\mathrm{d}\nu)$$ is an isometrical
isomorphism of the space $L_2(\mathbb{R}_+)$ onto $A^{(k)}$.
Clearly, see also Figure~\ref{fig_operatorsIII}, operators $R_{k}$
and $R_{k}^{*}$ provide the following decompositions of the
projection $P^{(k)}$ and of the identity operator on
$L_{2}(\mathbb{R}_+)$, i.e.,
\begin{align*}
R_{k} R_{k}^{*} & = I_{\phantom{k}} : L_{2}(\mathbb{R}_+) \to
L_{2}(\mathbb{R}_+),
\\ R_{k}^{*} R_{k} & = P^{(k)} : L_{2}(G,\mathrm{d}\nu) \to
A^{(k)}.
\end{align*}By a direct computation we have the explicit forms of
both operators providing Bargmann-type transforms in our
situation.

\begin{thm}\label{IIR*}
The isometrical isomorphism $R_{k}^{*}=U^{*} Q_k:
L_{2}(\mathbb{R}_+) \to A^{(k)}$ is given by
\begin{equation}\label{R*}
\left(R_{k}^{*}f\right)(\zeta) = \sqrt{2}\,v \int_{\mathbb{R}_+}
f(\xi)\ell_k(2\xi v)\,\mathrm{e}^{2\pi \mathrm{i}\xi
u}\,\sqrt{\xi}\,\mathrm{d}\xi,
\end{equation}with $\zeta=(u,v)\in G$.
The inverse isomorphism $R_{k}=Q_k^{*} U: A^{(k)} \to
L_{2}(\mathbb{R}_+)$ has the following form
\begin{equation}\label{RF}
(R_k F)(\xi) = \chi_+(\xi)
\sqrt{2\xi}\int_{\mathbb{R}\times\mathbb{R}_+} F(u,v)
\ell_k(2v\xi)\,\mathrm{e}^{-2\pi \mathrm{i}\xi
u}\,\frac{\mathrm{d}u\mathrm{d}v}{v}.
\end{equation}
\end{thm}

An interesting question how these operators are related to induced
representations of the affine group $G$ is solved in recent
paper~\cite{EHK}.

\section{Calder\'on-Toeplitz operators}\label{section3}

Wavelet transforms, including the ones coming from the affine
group $G$, are the building block of localization operators, see
book~\cite{wongbook} for further details on wavelet transforms and
localization operators. The representation of wavelet subspaces
summarized in previous section is especially important in the
study of Toeplitz-type operators related to wavelets which symbols
depend only on vertical variable $v=\Im \zeta$ in the upper
half-plane $\Pi$ of the complex plane $\mathbb{C}$,
see~\cite{hutnik2} for the general setting. This "restriction" to
imaginary part of a complex number is due to the decomposition
scheme we have just used, but on the other hand, it allows to
investigate properties of Toeplitz-type operators in a very
elegant way. Moreover, it gives rise to commutative algebras of
these operators in both cases of bounded and also unbounded
symbols which will be of further interest elsewhere.

For a given bounded function $a$ on $G$ define the Toeplitz-type
operator $T_a^{(k)}: A^{(k)}\to A^{(k)}$ with symbol $a$ as
$$T_{a}^{(k)} = P^{(k)}M_{a},$$ where $M_a$ is the operator
of pointwise multiplication by $a$ on $L_2(G,\mathrm{d}\nu)$ and
$P^{(k)}$ is the orthogonal projection from $L_2(G,\mathrm{d}\nu)$
onto $A^{(k)}$.

\begin{rem}\rm
In what follows we always consider the operators $T_a^{(k)}$
acting on wavelet subspaces $A^{(k)}$ although we may also define
the operators $\bar{T}_a^{(k)}$ acting on $\bar{A}^{(k)}$ and
given by $\bar{T}_{a}^{(k)} = \bar{P}^{(k)}M_{a}$. It is worth
noting that in this case of many wavelet subspaces (parameterized
by $k$) other Toeplitz- and Hankel-type (or, Ha-plitz in the
terminology of paper~\cite{nikolski}) operators may be defined as
follows
\begin{align*}
T_a^{(k,l)} & = P^{(k)}M_a P^{(l)}, \\
h_a^{(k,l)} & = \bar{P}^{(k)}M_a P^{(l)}, \\
H_a^{(k,l)} & = \left(I-\sum_{j=0}^{k}P^{(j)}\right)M_a P^{(l)},
\end{align*}see the works~\cite{jiang-peng} and~\cite{peng}.
\end{rem}

There exists an alternative way how to get the Calder\'on-Toeplitz
operators. In fact, $T_a^{(k)}$ may be viewed as operators acting
on $L_2(\mathbb{R})$ defined by the formula
$$\left\langle T_a^{(k)}f, g \right\rangle = \int_G a(\zeta) \left\langle f, \rho_\zeta\psi^{(k)}\right\rangle
\left\langle \rho_\zeta\psi^{(k)},
g\right\rangle\,\mathrm{d}\nu(\zeta), \quad f, g\in
L_2(\mathbb{R}),$$ interpreted in a weak sense. In this case the
identity $T_1^{(k)} f=f$ is known as the \textit{Calder\'on
reproducing formula}, cf.~\cite{calderon}, usually being used to
define classes of Hilbert spaces with reproducing kernels.
Therefore these operators are known as the
\textit{Calder\'on-Toeplitz operators}, cf.~\cite{nowak3}, and
were introduced by Richard Rochberg in~\cite{rochberg-CT} as a
wavelet counterpart of~Toeplitz operators defined on Hilbert
spaces of~holomorphic functions. Also, they are an effective
time-frequency localization tool in the context of wavelet
analysis, see~\cite{daubechies}, which provides ways of analyzing
signals by describing their frequency content as it varies over
time, and therefore they are a natural counterpart to the
intensively studied localization operators in time-frequency
analysis, see e.g.~\cite{CG}. From it follows that the
Calder\'on-Toeplitz operator $T_a^{(k)}$ may be viewed as a
filtered version of a signal $f$ with a symbol $a$ being
considered as a time-varying filter emphasizing or eliminating
some kind of information contained in time-scale content on level
$k$. For further information and results for Calder\'on-Toeplitz
operators we refer to papers of Nowak~\cite{nowak3},
\cite{nowak2}, \cite{nowak4}, and of Rochberg~\cite{rochberg-CT},
\cite{rochberg2}, \cite{rochberg3}.

Now we will demonstrate the usefulness of Bargmann-type transform
$R_k$ and its inverse $R_k^*$ (given by~$\mathrm{(\ref{RF})}$
and~$\mathrm{(\ref{R*})}$) under which we may study the unitary
equivalent images $R_k T_a^{(k)} R_k^*$ of Calder\'on-Toeplitz
operators $T_a^{(k)}$, see also the paper~\cite{huthut} for a
slightly different approach, but more general results obtained
therein. For our investigations the following result is very
important because in the case of symbols depending on vertical
coordinate of $G$ it enables us to reduce the Calder\'on-Toeplitz
operator to a certain multiplication operator. For the sake of
completeness we give its short proof here.

\begin{thm}\label{CTO1}
Let $(u,v)\in G$. If a measurable symbol $a=a(v)$ does not depend
on $u$, then the Calder\'on-Toeplitz operator $T_{a}^{(k)}$ acting
on $A^{(k)}$ is unitarily equivalent to the multiplication
operator $\mathfrak{A}_a^{(k)} = \gamma_{a,k} I$ acting on
$L_2(\mathbb{R}_+)$, where the function $\gamma_{a,k}:
\mathbb{R}_+\to\mathbb{C}$ is given by
\begin{equation}\label{gamma1}
\gamma_{a,k}(\xi) =
\int_{\mathbb{R}_{+}}a\left(\frac{v}{2\xi}\right)\ell_k^2(v)\,\mathrm{d}v,
\quad \xi\in\mathbb{R}_+.
\end{equation}
\end{thm}

\proof From construction of operators presented in the previous
section we directly have the following sequence of operator
equalities
\begin{align*}
\mathfrak{A}_a^{(k)} & = R_k T^{(k)}_{a} R_k^* = R_k P^{(k)}M_{a}
P^{(k)} R_k^* = R_k (R_k^* R_k) a(v) (R_k^* R_k) R_k^* \\ & = (R_k
R_k^*) R_k a(v) R_k^*(R_k R_k^*) = R_k a(v) R_k^* \\ & = Q_k^* U_2
U_1 a(v) U_1^{-1} U_2^{-1} Q_k = Q_k^* U_2 a(v) U_2^{-1} Q_k.
\end{align*}Since for a function $F\in L_2(\mathbb{R}, \mathrm{d}x)\otimes L_2(\mathbb{R}_+, \mathrm{d}y)$
holds $$\Bigl(U_2\,a(v) U_2^{-1}F\Bigr)(x,y) = U_2\left(a(v)
\sqrt{2|\omega|}\,v F(\omega, 2|\omega|v)\right) =
a\left(\frac{y}{2|x|}\right) F(x,y),$$ then
\begin{align*}
\left(\mathfrak{A}_a^{(k)}f\right)(\xi) & = \Bigl(Q_k^* U_2\,a(v)
U_2^{-1} Q_k f\Bigr)(\xi) \\ & = Q_k^* \left[\chi_+(x)\,
a\left(\frac{y}{2|x|}\right) f(x) \ell_k(y)\right](\xi) \\ & =
f(\xi)\, \chi_+(\xi) \int_{\mathbb{R}_+}
a\left(\frac{t}{2|\xi|}\right)\,\ell_k^2(t)\,\mathrm{d}t \\ & =
\gamma_{a,k}(\xi) f(\xi), \quad \xi\in\mathbb{R}_+,
\end{align*}which completes the proof.
\qed

We may observe that the function $\gamma_{a,k}$ is constructed by
putting a multiplier in admissibility
condition~(\ref{admissibility}) with respect to wavelet
$\psi^{(k)}$. As we will see later the function $\gamma_{a,k}$ is
responsible for many properties of the corresponding
Calder\'on-Toeplitz operator $T_a^{(k)}$ with a symbol $a=a(v)$
(bounded, but also unbounded one), and it shed a new light upon
the investigation of properties of the corresponding
Calder\'on-Toeplitz operator.

Let us mention that for general symbols $a=a(u,v)$ on $G$ the
Calder\'on-Toeplitz operator $T_a^{(k)}$ is no longer unitarily
equivalent to a multiplication operator. In fact, the operator
$R_k T_a^{(k)} R_k^*$ may have a more complicated structure: we
clarify this statement for a symbol depending only on the first
individual coordinate of $G$. This follows from our recent
paper~\cite{huthut} as a special case of Theorem~3.8 presented
therein, but for the sake of completeness and a little different
approach we present here the formal computations providing its
proof. In what follows put $\mathbb{R}_+^2:=
\mathbb{R}_+\times\mathbb{R}_+$.

\begin{thm}\label{CTO2} Let $(u,v)\in G$.
If a measurable function $b=b(u)$ does not depend on $v$, then the
Calder\'on-Toeplitz operator $T_{b}^{(k)}$ acting on $A^{(k)}$ is
unitarily equivalent to the operator $\mathfrak{B}_b^{(k)}$ acting
on $L_2(\mathbb{R}_+)$ given by
$$\left(\mathfrak{B}_b^{(k)}f\right)(\xi) = \int_{\mathbb{R}_+}
\mathcal{B}_k(\xi,t)\, \hat{b}(\xi-t) f(t)\,\mathrm{d}t, \quad
\xi\in\mathbb{R}_+,
$$ where the function $\mathcal{B}_k:
\mathbb{R}_+^2\to\mathbb{C}$ has the form
\begin{equation}\label{function B_k}
\mathcal{B}_k(\xi,t) = \frac{2\sqrt{t\xi}}{t+\xi}\,
P_k\left(\frac{8t\xi}{(t+\xi)^2}-1\right)\end{equation} with
$$P_n(x) = \frac{1}{2^n n!}\, \frac{\mathrm{d}^n}{\mathrm{d}x^n}
(x^2-1)^n$$ being the Legendre polynomial of degree
$n\in\mathbb{Z}_+$ for $x\in [-1,1]$.
\end{thm}

\proof Similarly as in the proof of Theorem~\ref{CTO1}, for a
function $f\in L_2(\mathbb{R}_+)$ we have
\begin{align*}
\mathfrak{B}_b^{(k)} f & = R_k T^{(k)}_{b} R_k^* f = 
R_k b(u) R_k^* f \\ & = Q_k^* U_2 (\mathcal{F}\otimes I) b(u)
(\mathcal{F}^{-1}\otimes I) U_2^{-1} Q_k f \\
& = 
Q_k^* U_2 (\mathcal{F}\otimes I) b(u) (\mathcal{F}^{-1}\otimes I)
\left(\chi_+(\omega)\sqrt{2|\omega|}\,v
f(\omega)\ell_k(2|\omega|v)\right).\end{align*}Using the
convolution theorem for Fourier transform we get
\begin{align*}
\left(\mathfrak{B}_b^{(k)} f\right)(\xi) & = Q_k^* U_2\left(
\int_{\mathbb{R}} \chi_+(t)\sqrt{2|t|}\,v\, \hat{b}(\omega-t)
f(t)\ell_k(2|t|v)\,\mathrm{d}t\right)(\xi)
\\ & = Q_k^* \left(\int_{\mathbb{R}_+}
\sqrt{\frac{t}{|x|}}\,\hat{b}(x-t)
f(t)\ell_k\left(\frac{ty}{|x|}\right)\,\mathrm{d}t\right)(\xi) \\
& = \chi_+(\xi)\int_{\mathbb{R}_+}
\sqrt{\frac{t}{|\xi|}}\,\hat{b}(\xi-t) f(t)
\left(\int_{\mathbb{R}_+} \ell_k(y)
\ell_k\left(\frac{ty}{|\xi|}\right)\,\mathrm{d}y\right)\,\mathrm{d}t
\\ & = \chi_+(\xi)\int_{\mathbb{R}_+} \hat{b}(\xi-t) f(t)
\left(\sqrt{t|\xi|}\,\int_{\mathbb{R}_+}
\ell_k(t\tau)\ell_k(|\xi|\tau)\,\mathrm{d}\tau\right)\,\mathrm{d}t.
\end{align*}Putting $$\mathcal{B}_k(\xi,t) = \chi_+(\xi) \sqrt{t|\xi|}
\int_{\mathbb{R}_+}
\ell_k(t\tau)\ell_k(|\xi|\tau)\,\mathrm{d}\tau, \quad
t\in\mathbb{R}_+,$$ we have the desired result. It remains to show
that $\mathcal{B}_k$ has the explicit form given by the
formula~(\ref{function B_k}). To this end we use~\cite[formula
7.414.13]{GR} to obtain
\begin{align*}
\mathcal{B}_k(\xi,t) & = \sqrt{t\xi} \int_{\mathbb{R}_+}
\exp\left[-\tau\left(\frac{t+\xi}{2}\right)\right]L_k(t\tau)
L_k(\xi\tau)\,\mathrm{d}\tau \\ & =\frac{2\sqrt{t\xi}}{t+\xi}\,
P_k\left(\frac{8t\xi}{(t+\xi)^2}-1\right), \quad
(t,\xi)\in\mathbb{R}_+^2,\end{align*}and the proof is complete.
\qed

As can be seen, the results of Theorem~\ref{CTO1} and
Theorem~\ref{CTO2} involving Calder\'on-Toeplitz operators with
symbols depending on the individual coordinates of $G$ describe an
analogy between the Calder\'on-Toeplitz operators and the calculus
of pseudo-differential operators. We clarify this fact for the
case of symbol in the product form. For the more general result
see~\cite[Theorem 3.10]{huthut}.

\begin{thm}
Let $(u,v)\in G$. If $c(u,v)=a(v)b(u)$ is a measurable symbol,
then the Calder\'on-Toeplitz operator $T_{c}^{(k)}$ acting on
$A^{(k)}$ is unitarily equivalent to the pseudo-differential
operator $\mathfrak{C}_{\mathfrak{c}}^{(k)}$ acting on
$L_2(\mathbb{R}_+)$ given by the iterated integral
$$\left(\mathfrak{C}_{\mathfrak{c}}^{(k)} f\right)(\xi) = \int_{\mathbb{R}} du
\int_{\mathbb{R}_+} \mathfrak{c}_k(\xi,t,u) f(t)\,
\mathrm{e}^{-2\pi \mathrm{i}(\xi-t)u}\,\mathrm{d}t, \quad
\xi\in\mathbb{R}_+,
$$ where its compound (double) symbol $\mathfrak{c}_k: \mathbb{R}_+^2\times\mathbb{R}\to\mathbb{C}$
has the form
$$\mathfrak{c}_k(\xi,t,u) = 2\sqrt{t\xi}\ b(u) \int_{\mathbb{R}_+}
a(v) \ell_k(2v\xi) \ell_k(2vt)\,\mathrm{d}v.$$
\end{thm}


\begin{rem}\rm
For each $k\in\mathbb{Z}_+$ the function $\mathcal{C}_{a,k}:
\mathbb{R}_+^2 \to \mathbb{C}$ given in the form of improper
parameter-dependent integral
$$\mathcal{C}_{a,k}(\xi,t) = 2\sqrt{t\xi}\ \int_{\mathbb{R}_+}
a(v) \ell_k(2v\xi) \ell_k(2vt)\,\mathrm{d}v$$ provides an
"extension" of both functions $\gamma_{a,k}$ and $\mathcal{B}_k$.
Indeed, in the first case $\gamma_{a,k}$ is a restriction of
$\mathcal{C}_{a,k}$ to the diagonal, i.e.,
$\mathcal{C}_{a,k}(\xi,\xi)=\gamma_{a,k}(\xi)$ for each
$\xi\in\mathbb{R}_+$, and in the second case $\mathcal{B}_k$ is a
restriction of $\mathcal{C}_{a,k}$ to a constant symbol $a$, i.e.,
$\mathcal{C}_{1,k}(\xi,t)=\mathcal{B}_k(\xi,t)$ for each
$(\xi,t)\in\mathbb{R}_+^2$. Also it can be proved that for a
bounded function $a=a(v)$ the function $\mathcal{C}_{a,k}$ is
continuous and bounded on $\mathbb{R}_+^2$ for each
$k\in\mathbb{Z}_+$.
\end{rem}

In what follows we return back to the case of symbols depending on
$v\in\mathbb{R}_+$ which provides a number of results for
properties of Calder\'on-Toeplitz operators and their algebras as
direct corollaries of Theorem~\ref{CTO1}.

Given a linear subset $\mathcal{A}$ of $L_\infty(\mathbb{R}_+)$,
for $k\in\mathbb{Z}_+$ denote by $\mathcal{T}_k(\mathcal{A})$ the
$C^*$-algebra generated by all Calder\'on-Toeplitz operators
$T_a^{(k)}$ with symbols $a\in\mathcal{A}$ acting on the wavelet
subspace $A^{(k)}$. As a first useful algebra of symbols we
introduce the $C^*$-algebra $\mathcal{A}_\infty$ of all bounded
functions on $G$ depending only on $v = \Im\zeta$, $\zeta\in G$.
Then the following result is in the spirit of Vasilevski results,
see e.g.~\cite[Corollary 10.4.10]{vasilevskibook}, obtained for
Toeplitz operators on (weighted) Bergman spaces. Put
$C_b(\mathbb{R}_+):=C(\mathbb{R}_+)\cap L_\infty(\mathbb{R}_+)$.

\begin{cor}\rm\label{coralgebra}
Each $C^*$-algebra $\mathcal{T}_k(\mathcal{A}_\infty)$,
$k\in\mathbb{Z}_+$, is commutative and is isometrically imbedded
to the algebra $C_b(\mathbb{R}_+)$. The isomorphic imbedding
$$\tau_\infty^{(k)}:
\mathcal{T}_k(\mathcal{A}_\infty)\longrightarrow
C_b(\mathbb{R}_+)$$ is generated by the following mapping of
generators of $\mathcal{T}_k(\mathcal{A}_\infty)$
$$\tau_\infty^{(k)}: T_a^{(k)} \longmapsto \gamma_{a,k}(\xi),$$ where $a\in
\mathcal{A}_\infty$.
\end{cor}


Property to be unitarily equivalent to a multiplication operator
permits us to describe easily invariant subspaces of algebra
$\mathcal{T}_k(\mathcal{A}_\infty)$. Note that the following
result still holds for any $C^*$-algebra generated by bounded
Calder\'on-Toeplitz operators with unbounded (measurable) symbol
depending on the imaginary part of complex number. A more detailed
study of boundedness of Calder\'on-Toeplitz operators with
unbounded symbols is done in paper~\cite{hutnikII}.

\begin{cor}\rm
Each commutative $C^*$-algebra $\mathcal{T}_k(\mathcal{A}_\infty)$
is reducible. Every invariant subspace $\mathcal{S}_k$ of
$\mathcal{T}_k(\mathcal{A}_\infty)$ is defined by a measurable
subset $S_k\subset \mathbb{R}_+$ and has the form
$$\mathcal{S}_k=(R_k^*\chi_{S_k} I)L_2(\mathbb{R}_+).$$
\end{cor}

Reverting the statement of Theorem~\ref{CTO1} we get the following
spectral-type representation of a Calder\'on-Toeplitz operator.
Its proof goes directly from Theorem~\ref{CTO1} and
Theorem~\ref{IIR*}.

\begin{cor}\label{representationCTO}\rm
For $a\in\mathcal{A}_\infty$ the Calder\'on-Toeplitz operator
$T_{a}^{(k)}$ acting on $A^{(k)}$ admits the following
representation
\begin{equation}\label{repCTO}
\left(T_{a}^{(k)}F\right)(\zeta) =
\sqrt{2}\,v\int_{\mathbb{R}_+}\gamma_{a,k}(\xi)\ell_k(2v\xi)f(\xi)\,\mathrm{e}^{2\pi
\mathrm{i}u\xi}\,\sqrt{\xi}\,\mathrm{d}\xi,
\end{equation} where $\zeta=(u,v)\in G$ and $f(\xi)=(R_k F)(\xi)$.
\end{cor}

It may be observed that all the above stated constructions fit
perfectly to the general coherent states scheme summarized e.g.
in~\cite[Appendix A]{vasilevskibook} (i.e., wavelets = affine
coherent states). Following this scheme the next result gives the
form of the Wick symbol of Calder\'on-Toeplitz operator
$T_{a}^{(k)}$ depending on $v=\Im\zeta$. Note that writing the
Calder\'on-Toeplitz operator $T_{a}^{(k)}$ in terms of its Wick
symbol yields exactly the spectral-type
representation~(\ref{repCTO}). Recall that $\kappa_k$ is the
constant depending on $k$ given in~(\ref{lambda_k}).

\begin{cor}\label{thmwicksymbol}\rm
Let $\zeta=(u,v)\in G$. Given $a=a(v)\in\mathcal{A}_\infty$, the
Wick symbol $\widetilde{a}_k(\zeta)$ of the Calder\'on-Toeplitz
operator $T_{a}^{(k)}$ depends only on $v$ as well, and has the
form
\begin{equation}\label{wicksymbol}
\widetilde{a}_k(v) = 2\,\kappa_k^{-1} v^2\int_{\mathbb{R}_+}
\gamma_{a,k}(\xi)\ell_k^2(2v\xi)\,\xi \mathrm{d}\xi.
\end{equation}
The corresponding Wick function is given by the formula
$$\widetilde{a}_k(\zeta,\eta) = \frac{2tv}{K^{(k)}_{\zeta}(\eta)}
\int_{\mathbb{R}_+}\gamma_{a,k}(\xi)\ell_k(2v\xi)\ell_k(2t\xi)
\,\mathrm{e}^{-2\pi \mathrm{i}\xi(u-s)}\,\xi \mathrm{d}\xi, $$with
$\zeta=(u,v)$, $\eta=(s,t) \in G$.
\end{cor}

\begin{rem}\rm
Formula~(\ref{wicksymbol}) may be interpreted in the following
interesting way. For the Calder\'on-Toeplitz operator $T_a^{(k)}$
acting on $A^{(k)}$ with a symbol $a\in\mathcal{A}_\infty$
calculate the corresponding function $\gamma_{a,k}(\xi)$,
$\xi\in\mathbb{R}_+$. Let us introduce the multiplication operator
$$\left(M^{(k)}_a f\right)(x) = \kappa_k^{-1/2}a(x)f(x), \quad
f\in L_2(\mathbb{R}_+).$$ Take the function
$\left(M^{(k)}_{\textrm{Id}}\gamma_{a,k}\right)(\xi)$, where
$\textrm{Id}$ is the identity function, and consider the
Calder\'on-Toeplitz operator on $A^{(k)}$ with symbol
$\left(M^{(k)}_{\textrm{Id}}\gamma_{a,k}\right)(\xi)$. Then the
function
$$\left(M^{(k)}_{\textrm{Id}}\gamma_{\left(M^{(k)}_{\textrm{Id}}\gamma_{a,k}\right)(\xi),k}\right)(v)$$ is
nothing but the Wick symbol of the initial Calder\'on-Toeplitz
operator $T_a^{(k)}$, i.e.,
$$\widetilde{a}_k(v) = \widetilde{a}_k(\zeta, \zeta) =
\left(M^{(k)}_{\textrm{Id}}\gamma_{\left(M^{(k)}_{\textrm{Id}}\gamma_{a,k}\right)(\xi),k}\right)(v),
\quad \zeta=(u,v)\in G.$$
\end{rem}

Let us mention the star product $\star$ defining the composition
of two Wick symbols $\widetilde{a}_A$ and $\widetilde{a}_B$ of two
operators $A$ and $B$ as the Wick symbol of their composition
$AB$, i.e., $\widetilde{a}_A\star\widetilde{a}_B =
\widetilde{a}_{AB}$. The following result gives the formula for
the star product of two Calder\'on-Toeplitz operators in terms of
the corresponding function $\gamma$. Again, it is an immediate
consequence of Theorem~\ref{CTO1} and
Corollary~\ref{thmwicksymbol}.

\begin{cor}\rm
Let $(u,v)\in G$. Let $T_{a}^{(k)}$ and $T_{b}^{(k)}$ be two
Calder\'on-Toeplitz operators acting on $A^{(k)}$ with symbols
$a(v)$ and $b(v)$, and let $\widetilde{a}_k(v)$ and
$\widetilde{b}_k(v)$ be their Wick symbols, respectively. Then the
Wick symbol $\widetilde{c_k}$ of the composition
$T_{a}^{(k)}T_{b}^{(k)}$ is given by
$$\widetilde{c_k}(v) = \left(\widetilde{a}_k \star \widetilde{b}_k\right)(v)
= 2\,\kappa_k^{-1}v^2
\int_{\mathbb{R}_+}\gamma_{a,k}(\xi)\gamma_{b,k}(\xi)\ell_k^2(2v\xi)\,\xi
\mathrm{d}\xi.$$
\end{cor}

\section{Some properties of function $\gamma_{a,k}$}\label{section4}

In this section we investigate basic properties of function
$\gamma_{a,k}: \mathbb{R}_+\to\mathbb{C}$ related to
Calder\'on-Toeplitz operator $T_a^{(k)}$ with a symbol $a=a(v)$ in
detail. Firstly, let us summarize that for each $k\in\mathbb{Z}_+$
and each $a\in \mathfrak{A}_\infty$ we have $\gamma_{a,k}\in
L_\infty(\mathbb{R}_+)$ (for more details on this "boundedness
topic" and its consequences for Calder\'on-Toeplitz operators see
the paper~\cite{hutnikII}). Moreover, in such a case of bounded
symbol $a$ the function $\gamma_{a,k}(\xi)$ is also continuous in
each finite point $\xi\in\mathbb{R}_+$, and thus $\gamma_{a,k}\in
C_b(\mathbb{R}_+)$ -- this fact was already used in
Corollary~\ref{coralgebra}. The continuity, but also other
properties of $\gamma_{a,k}$, may be perhaps better seen if we use
some formulas for Laguerre polynomials, more
precisely~\cite[formula 8.976.3]{GR} and \cite[formula (5), p.
209]{Rainville} combined with the exact form of $L_n$ given
in~(\ref{defLaguerre}) enable to rewrite the function
$\gamma_{a,k}(\xi)$ in the form
$$\gamma_{a,k}(\xi) = \sum_{i=0}^k \sum_{j=0}^{2k} \sum_{r=0}^{j}
c(k,i,j,r)(1-4\xi)^{2i-j} (4\xi)^{j+1} \int_{\mathbb{R}_+} a(v)v^r
\mathrm{e}^{-2v\xi}\,\mathrm{d}v,$$ where
$$c(k,i,j,r)=\frac{1}{2^{2k+1}} {2k-2i \choose k-i} {2i \choose j}
{j\choose r} {2i\choose i}\frac{(-1)^r }{r!}.$$ Here again might
be seen the close connection with the classical Toeplitz operators
acting on weighted Bergman spaces (the so called parabolic case),
because the last formula resembles the function $\gamma_{a,
\lambda}$ obtained in that case by Vasilevski and his
collaborators, see~\cite[formula (2.6)]{GKV}, or~\cite[formula
(10.4.4), p. 254]{vasilevskibook}.

Now, we describe the behavior of higher order derivatives of
function $$\gamma_{a,k}(\xi) = 2\xi \int_{\mathbb{R}_+} a(v)
\ell_k^2(2v\xi)\,\mathrm{d}v$$ as $\xi\to+\infty$. Trivially, for
a constant symbol $a$ on $\mathbb{R}_+$ the function
$\gamma_{a,k}$ is constant on $\mathbb{R}_+$ for each
$k\in\mathbb{Z}_+$, and therefore all its higher order derivatives
are zero. Let us consider a non-constant bounded function $a=a(v)$
on $\mathbb{R}_+$. Differentiating $n$-times (with $n\geq 1$) we
obtain
\begin{equation}\label{n-thderivativegamma}
\frac{\mathrm{d}^n\gamma_{a,k}(\xi)}{\mathrm{d}\xi^n} =
2n\int_{\mathbb{R}_+}
a(v)\frac{\mathrm{d}^{n-1}}{\mathrm{d}\xi^{n-1}}\ell_k^2(2v\xi)\,\mathrm{d}v
+ 2\xi \int_{\mathbb{R}_+} a(v)
\frac{\mathrm{d}^n}{\mathrm{d}\xi^n}\ell_k^2(2v\xi)\,\mathrm{d}v,\end{equation}
and thus,
$$\left|\frac{\mathrm{d}^n\gamma_{a,k}(\xi)}{\mathrm{d}\xi^n}\right|
\leq C \left(\frac{n}{\xi}
I_{k}^{(n-1)}(\xi)+I_k^{(n)}(\xi)\right),$$ where
$$I_k^{(m)}(\xi) := \int_{\mathbb{R}_+}
\left|\frac{\mathrm{d}^{m}}{\mathrm{d}\xi^{m}}\ell_k^2(2v\xi)\right|\,2\xi\mathrm{d}v,
\quad m\in\mathbb{Z}_+.$$ Using the
formula~(\ref{n-thderivativeell}) and function $\Lambda$ from
Appendix, the function $I_k^{(m)}$ may be rewritten as follows
\begin{align*}
I_k^{(m)}(\xi) & = \frac{1}{\xi^m} \sum_{i=0}^m \sum_{j=0}^i {m
\choose i} {i \choose j} \int_{\mathbb{R}_+} (2v\xi)^m
\mathrm{e}^{-2v\xi} \left|L_{k-i+j}^{(i-j)}(2v\xi)
L_{k-j}^{(j)}(2v\xi)\right| 2\xi\mathrm{d}v \\ & = \frac{1}{\xi^m}
\sum_{i=0}^m \sum_{j=0}^i {m \choose i} {i \choose j}
\int_{\mathbb{R}_+}
\Lambda_{m,k-i+j,k-j}^{(i-j,j)}(x)\,\mathrm{d}x.
\end{align*}Since the last integral is finite, see the formula~(\ref{upperbound2})
from Appendix, then $I_k^{(m)}(\xi)\to 0$ as $\xi\to +\infty$ for
each $m\geq 1$ and each $k\in\mathbb{Z}_+$. Clearly, in the case
$m=0$ we have $$I_k^{(0)}(\xi)= \int_{\mathbb{R}_+}
\Lambda_{0,k,k}^{(0,0)}(x)\,\mathrm{d}x = \int_{\mathbb{R}_+}
\ell_k^2(x)\,\mathrm{d}x = 1$$ for each $k\in\mathbb{Z}_+$, and
thus $\xi^{-1}I_k^{(0)}(\xi)\to 0$ as $\xi\to+\infty$. So, for
$a=a(v)\in L_\infty(\mathbb{R}_+)$ we get that for each
$n=1,2,\dots$ holds
\begin{equation}\label{derivativegamma}
\lim_{\xi\to+\infty}\frac{\mathrm{d}^n\gamma_{a,k}(\xi)}{\mathrm{d}\xi^n}
= 0\end{equation} for each $k\in\mathbb{Z}_+$. One can observe
that for a non-constant bounded symbol $a$ on $\mathbb{R}_+$ each
derivative of $\gamma_{a,k}(\xi)$ \textit{is unbounded} near the
point $\xi=0$ for each $k\in\mathbb{Z}_+$.

In the following theorem we show that the above observation for
behavior of higher order derivatives of $\gamma_{a,k}$ holds also
if we replace a bounded function $a=a(v)$ by an $L_1$-integrable
function on $\mathbb{R}_+$. The proof requires a number of
estimates which are for better readability given in Appendix.

\begin{thm}
Let $(u,v)\in G$ and $a=a(v)\in L_1(\mathbb{R}_+)$ be such that
$\gamma_{a,k}(\xi)\in L_\infty(\mathbb{R}_+)$. Then for each
$n=1,2,\dots$ the equation~$\mathrm{(\ref{derivativegamma})}$
holds for each $k\in\mathbb{Z}_+$.
\end{thm}

\proof Let $n\geq 1$. Then differentiating $n$-times
yields~(\ref{n-thderivativegamma}), and using the
formula~(\ref{n-thderivativeell}) from Appendix we have
\begin{align}\label{dngamma}
\frac{\mathrm{d}^n\gamma_{a,k}(\xi)}{\mathrm{d}\xi^n} & =
n(-1)^{n-1}\sum_{i=0}^{n-1} \sum_{j=0}^i {n-1 \choose i} {i
\choose j} \frac{I^{(n-1)}_{k,i,j}(\xi)}{\xi} \nonumber \\ &
\phantom{=} + (-1)^{n} \sum_{i=0}^n \sum_{j=0}^i {n \choose i} {i
\choose j} I^{(n)}_{k,i,j}(\xi),\end{align}where
$$I^{(m)}_{k,i,j}(\xi) := \frac{1}{\xi^{m}} \int_{\mathbb{R}_+}
a(v)(2v\xi)^m \mathrm{e}^{-2v\xi} L_{k-i+j}^{(i-j)}(2v\xi)
L_{k-j}^{(j)}(2v\xi)\,2\xi\mathrm{d}v, \quad m\in\mathbb{Z}_+.$$
Observe that for $n=1$ the first term in~(\ref{dngamma}) is equal
to $2\xi^{-1}\gamma_{a,k}(\xi)$ which tends to 0 as $\xi\to
+\infty$. We now show that also all the above integrals
$I_{k,i,j}^{(m)}(\xi)$ for $m\geq 1$ tend to $0$ as
$\xi\to+\infty$.

Let $m\geq 1$. For a sufficiently small $\delta>0$ consider the
integral
\begin{align*}
I_{k,i,j}^{(m)}(\xi) & = \frac{1}{\xi^m}\Biggl(\int_{0}^\delta
a(v)(2v\xi)^m \mathrm{e}^{-2v\xi} L_{k-i+j}^{(i-j)}(2v\xi)
L_{k-j}^{(j)}(2v\xi)\,2\xi\mathrm{d}v \\ & + \int_\delta^\infty
a(v)(2v\xi)^m \mathrm{e}^{-2v\xi}
L_{k-i+j}^{(i-j)}(2v\xi) L_{k-j}^{(j)}(2v\xi)\,2\xi\mathrm{d}v\Biggr) \\
& := I^{(m)}_{k,i,j,1}(\xi) + I^{(m)}_{k,i,j,2}(\xi).
\end{align*}Using~(\ref{upperbound2}) from Appendix we have
\begin{align*}
\left|I^{(m)}_{k,i,j,1}(\xi)\right| & \leq \sup_{v\in (0,\delta)}
|a(v)|\,\frac{1}{\xi^m} \int_{0}^\delta (2v\xi)^m
\mathrm{e}^{-2v\xi} \left|L_{k-i+j}^{(i-j)}(2v\xi)
L_{k-j}^{(j)}(2v\xi)\right|\,2\xi\mathrm{d}v \\ & = \sup_{v\in
(0,\delta)} |a(v)|\,\frac{1}{\xi^m}
\int_{0}^{2\delta\xi} \Lambda_{m,k-i+j,k-j}^{(i-j,j)}(x)\,\mathrm{d}x \\
& \leq \sup_{v\in (0,\delta)}
|a(v)|\,\frac{1}{\xi^m}\int_{\mathbb{R}_+}
\Lambda_{m,k-i+j,k-j}^{(i-j,j)}(x)\,\mathrm{d}x \\ & \leq
\sup_{v\in (0,\delta)} |a(v)|
\,\frac{\textrm{const}_{m,k-i+j,k-j}^{(i-j,j)}}{\xi^m}.
\end{align*}Thus for a sufficiently small $\delta$ and sufficiently large $\xi$ we
have $$\left|I^{(m)}_{k,i,j,1}(\xi)\right|<\varepsilon.$$ To
estimate $I^{(m)}_{k,i,j,2}$ we use the
formula~(\ref{upperbound1}) from Appendix to get
\begin{align*}
\left|I^{(m)}_{k,i,j,2}(\xi)\right| & \leq
\frac{2}{\xi^{m-1}}\int_\delta^\infty |a(v)|\,
\Lambda_{m,k-i+j,k-j}^{(i-j,j)}(2v\xi) \mathrm{d}v \\ & \leq
\sum_{r=0}^{k-i+j} \sum_{s=0}^{k-j}
\frac{(i-j+1)_{k-i+j-r}}{(k-i+j-r)!\,r!}
\frac{(j+1)_{k-j-s}}{(k-j-s)!\,s!}
\\ & \phantom{=} \times \frac{2}{\xi^{m-1}}\int_\delta^\infty |a(v)| (2v\xi)^{m+r+s}
\mathrm{e}^{-2v\xi}\,\mathrm{d}v
\\ & \leq \sum_{r=0}^{k-i+j} \sum_{s=0}^{k-j}
\frac{(i-j+1)_{k-i+j-r}}{(k-i+j-r)!\,r!}
\frac{(j+1)_{k-j-s}}{(k-j-s)!\,s!} \,
\\ & \phantom{=} \times \frac{2\mathrm{e}^{-\delta\xi}}{\xi^{m-1}}\int_\delta^\infty
|a(v)| (2v\xi)^{m+r+s} \mathrm{e}^{-v\xi}\,\mathrm{d}v.
\end{align*}Using~(\ref{estimation}) from Appendix we finally have
\begin{align*}
\left|I^{(m)}_{k,i,j,2}(\xi)\right| & \leq \sum_{r=0}^{k-i+j}
\sum_{s=0}^{k-j} \frac{(i-j+1)_{k-i+j-r}}{(k-i+j-r)!\,r!}
\frac{(j+1)_{k-j-s}}{(k-j-s)!\,s!} \\ & \phantom{=} \times
\left(\frac{2(m+r+s)}{\mathrm{e}}\right)^{m+r+s}\,
 \frac{2 \mathrm{e}^{-\delta\xi}}{\xi^{m-1}}\,\|a\|_{L_1(\mathbb{R}_+)},
\end{align*}which yields that for a sufficiently small $\delta$
and sufficiently large $\xi$ we have
$$\left|I^{(m)}_{k,j,i,2}(\xi)\right|<\varepsilon,$$ and therefore
$\lim\limits_{\xi\to+\infty} I^{(m)}_{k,i,j}(\xi)=0$ for all
$m\geq 1$. Finally, from the above also follows that all the
integrals $$J_{k,i,j}^{(m)}(\xi) :=
\frac{I_{k,i,j}^{(m)}(\xi)}{\xi}\to 0 \,\,\,\,\textrm{as}\,\,\,\,
\xi\to+\infty$$ for all $m\geq 1$, which finishes the proof. \qed

As a consequence we immediately have the following corollary.

\begin{cor}\rm
If $a=a(v)\in L_1(\mathbb{R}_+)\cup L_\infty(\mathbb{R}_+)$ such
that $\gamma_{a,k}\in L_\infty(\mathbb{R}_+)$, then for each
$n=1,2,\dots$ the equation~(\ref{derivativegamma}) holds for each
$k\in\mathbb{Z}_+$.
\end{cor}

We have already mentioned that for a bounded symbol $a$ the
function $\gamma_{a,k}$ is continuous on $\mathbb{R}_+$ for each
$k\in\mathbb{Z}_+$. Now we are interested in sufficient conditions
for its continuity on the whole $\overline{\mathbb{R}}_+ =
[0,+\infty]$ which will be useful when studying certain algebras
of operators. First observe that for a ,,very large $\xi$'' the
function $\ell^2_k(2v\xi)$ has a very sharp maximum at the point
$v=0$ and thus the major contribution to the integral for
$\xi\to+\infty$ is determined by values of $a(v)$ at a
neighborhood of the point $0$. On the other hand, the major
contribution for a ,,very small $\xi$'' is determined by values of
$a(v)$ at a neighborhood of $+\infty$. In particular, directly
from~(\ref{gamma1}) we deduce that if $a(v)$ has limits at the
points $0$ and $+\infty$, then
\begin{align*}
\lim_{\xi\to+\infty} \gamma_{a,k}(\xi) & = \lim_{v\to 0} a(v), \\
\lim_{\xi\to 0} \gamma_{a,k}(\xi) & = \lim_{v\to +\infty} a(v).
\end{align*}Now we do it precisely in the following theorem.

\begin{thm}\label{thmlimits}
Let $(u,v)\in G$. If $a=a(v)\in L_\infty(\mathbb{R}_+)$ and the
following limits exist
\begin{equation}\label{lima(v)}
\lim_{v\to 0} a(v)=a_0, \quad \lim_{v\to+\infty} a(v)=a_\infty,
\end{equation} then $\gamma_{a,k}(\xi)\in C[0,+\infty]$, and
$$\gamma_{a,k}(+\infty)=a_0, \quad \gamma_{a,k}(0)=a_\infty$$ for each
$k\in\mathbb{Z}_+$.
\end{thm}

\proof Let $\xi\to +\infty$. For a sufficiently small $\delta>0$
represent the function $\gamma_{a,k}(\xi)$ as follows
$$\gamma_{a,k}(\xi) = 2\xi\left(\int_{0}^\delta
a(v)\ell_k^2(2v\xi)\,\mathrm{d}v + \int_\delta^\infty
a(v)\ell_k^2(2v\xi)\,\mathrm{d}v\right) := I_1(\xi)+ I_2(\xi).
$$ Consider the integral $I_1$ in the form
$$I_1(\xi) = 2\xi \left(\int_{0}^\delta
a_0\ell_k^2(2v\xi)\,\mathrm{d}v + \int_{0}^\delta
[a(v)-a_0]\ell_k^2(2v\xi)\,\mathrm{d}v\right) :=
I_{1,1}(\xi)+I_{1,2}(\xi).$$ Then by Lemma~\ref{lemmaell} from
Appendix holds
$$I_{1,1}(\xi) = a_0 \int_{0}^{\delta} \ell_k^2(2v\xi)\,2\xi \mathrm{d}v =
a_0\int_{0}^{2\delta\xi} \ell_k^2(x)\,\mathrm{d}x =
a_0+\Phi(\delta,\xi),$$ where for a sufficiently large $\xi$ we
have
$$|\Phi(\delta,\xi)|<\frac{\varepsilon}{3},$$ and thus $|I_{1,1}(\xi)-a_0|<\frac{\varepsilon}{3}$.
Also, $$|I_{1,2}(\xi)| \leq \sup_{v\in (0,\delta)} |a(v)-a_0|
\int_{0}^{\delta} \ell_k^2(2v\xi)\,2\xi \mathrm{d}v = \sup_{v\in
(0,\delta)} |a(v)-a_0| \int_{0}^{2\delta\xi}
\ell_k^2(x)\,\mathrm{d}x,$$ which again means that for an
appropriate choice of $\delta$ and a sufficiently large $\xi$ we
have
$$|I_{1,2}(\xi)|<\frac{\varepsilon}{3}.$$
For the last integral $I_2$ we use the H\"older integral
inequality and Lemma~\ref{lemmaell} from Appendix to get
\begin{align*}
|I_2(\xi)| & \leq 2\xi \int_\delta^\infty |a(v)|
\ell_k^2(2v\xi)\,\mathrm{d}v \leq \|a\|_{L_\infty(\mathbb{R}_+)}
\int_\delta^\infty \ell_k^2(2v\xi)\,2\xi \mathrm{d}v \\ & =
\|a\|_{L_\infty(\mathbb{R}_+)} \left(1-\int_0^{2\delta\xi}
\ell_k^2(x)\,\mathrm{d}x\right) = \|a\|_{L_\infty(\mathbb{R}_+)}
N_{2k}(2\delta\xi)\,\mathrm{e}^{-2\delta \xi},
\end{align*}i.e., for a
sufficiently small $\delta$ and sufficiently large $\xi$ we get
$$|I_2(\xi)|<\frac{\varepsilon}{3}.$$ Summarizing the above we have that for
any $\varepsilon>0$ and an appropriate $\delta>0$ there is
$\xi_0>0$ such that for each $\xi\geq \xi_0$ one has
$$|\gamma_{a,k}(\xi)-a_0|<\varepsilon,$$ which proves that
$\lim\limits_{\xi\to+\infty} \gamma_{a,k}(\xi)=a_0$. The second
limit may be proved similarly. \qed

\begin{rem}\rm
The result of Theorem~\ref{thmlimits}, claiming that the limits at
infinity and at zero of the function $\gamma_{a,k}$ \textit{are
independent of parameter $k$}, is rather surprising in this case.
Since the wavelet transforms with Laguerre functions of order $k$
live on the true polyanalytic Bergman space of order $k$, then
from Theorem~\ref{thmlimits} follows that all the true
polyanalytic Bergman spaces have asymptotically the "same
behavior". The consequences of this observation in various areas
of mathematics and physics are not known for us at this moment.
\end{rem}

\begin{rem}\rm
Important assumption of Theorem~\ref{thmlimits} is boundedness of
symbol $a(v)$. Moreover, as the following example shows, the
condition~(\ref{lima(v)}) is \textit{not necessary} even for
bounded symbols: if $a(v)=\sin v$, $v\in\mathbb{R}_+$, then
$$\gamma_{a,1}(\xi) = 2\xi \int_{\mathbb{R}_+} \sin v\,
\mathrm{e}^{-2v\xi}(1-2v\xi)^2\,\mathrm{d}v =
\frac{2\xi\,(1-16\xi^2+48\xi^4)}{(1+4\xi^2)^3}, \quad
\xi\in\mathbb{R}_+,$$ which yields
$$\lim_{\xi\to +\infty} \gamma_{a,1}(\xi) = \lim_{\xi\to 0}
\gamma_{a,1}(\xi) = 0.$$
\end{rem}

On the other hand, in the following example we give an unbounded
symbol $a(v)$ such that the corresponding function
$\gamma_{\cdot}$ \textit{is continuous} on $[0,+\infty]$.

\begin{ex}\rm\label{exmunbounded}
The symbol $$a(v)=\frac{1}{\sqrt{v}}\sin\frac{1}{v}, \quad
v\in\mathbb{R}_+,$$ is clearly unbounded on $\mathbb{R}_+$. For
$k=1$ (we omit the details of lengthy computation) the
corresponding function $\gamma_{a,1}$ has the form
$$\gamma_{a,1}(\xi) = \frac{\sqrt{2\pi}}{4}\,\mathrm{e}^{-2\sqrt{\xi}}
\left[\left(2\sqrt{\xi}-8\xi\right)\frac{\cos
2\sqrt{\xi}}{2\sqrt{\xi}}+\left(3-2\sqrt{\xi}\right)\frac{\sin
2\sqrt{\xi}}{2\sqrt{\xi}}\right], $$ for $\xi\in\mathbb{R}_+$.
Then a direct computation yields $$\lim_{\xi\to 0}
\gamma_{a,1}(\xi) = \sqrt{2\pi} \,\,\,\,\,\textrm{and} \,\,\,\,\,
\lim_{\xi\to +\infty} \gamma_{a,1}(\xi) = 0,$$ which also means
that the function $\gamma_{a,1}$ is bounded on $\mathbb{R}_+$. For
$k\geq 2$ it is difficult to get an explicit form of
$\gamma_{a,k}$, so in~\cite{hutnikII} we present a different
approach to verify boundedness of $\gamma_{a,k}$, and
consequently, boundedness of the corresponding Calder\'on-Toeplitz
operator $T_a^{(k)}$ for each $k\in\mathbb{Z}_+$.
\end{ex}

\begin{ex}\rm\label{exmunbounded2}
The oscillating symbol $a(v)=\mathrm{e}^{2v \mathrm{i}}$ is
continuous at the point $v=0$. Therefore,
$\gamma_{a,k}(+\infty)=a(0)=1$ and it is sufficient to investigate
the behavior of $\gamma_{a,k}(\xi)$ as $\xi\to 0$. The explicit
form of the corresponding function $\gamma_{a,k}$ is as follows
\begin{align*}
\gamma_{a,k}(\xi) & = 2\xi \int_{\mathbb{R}_+} \mathrm{e}^{-2v
(\xi -\mathrm{i})} L_k^2(2v\xi)\,\mathrm{d}v =
\frac{\xi}{\xi-\mathrm{i}} \int_{\mathbb{R}_+} \mathrm{e}^{-t}
L_k^2\left(\frac{t\xi}{\xi-\mathrm{i}}\right)\,\mathrm{d}t \\ & =
\frac{(-1)^k}{(\xi-\mathrm{i})^{2k+1}}\sum_{j=0}^{k}
(-1)^{j}\left[{k\choose j}\right]^2 \xi^{2j+1}, \quad
\xi\in\mathbb{R}_+,
\end{align*} where~\cite[formula 7.414.2]{GR} has been used. Clearly, $\lim\limits_{\xi\to 0}
\gamma_{a,k}(\xi) = 0$, thus $\gamma_{a,k}(\xi)\in C[0,+\infty]$
for each $k\in\mathbb{Z}_+$.
\end{ex}

\section{Isomorphism between the Calder\'on-Toeplitz operator algebra and functional algebra}
\label{section5} 

In what follows we use the result of Theorem~\ref{thmlimits} to
study certain algebras of Calder\'on-Toeplitz operators. Therefore
we denote by $L^{\{0,+\infty\}}_\infty(\mathbb{R}_+)$ the
$C^*$-subalgebra of $L_\infty(\mathbb{R}_+)$ which consists of all
functions having limits at the points $0$ and $+\infty$.

\begin{rem}\rm
Considering a non-negative real-valued symbol $a$ from the algebra
$L^{\{0,+\infty\}}_\infty(\mathbb{R}_+)$ with a non-zero value
$a_0$, we get the following interesting property of function
$\gamma_{a,k}: \mathbb{R}_+\to \mathbb{R}_+$ for each
$k\in\mathbb{Z}_+$ in the form
$$\lim_{\xi\to\infty} \frac{\gamma_{a,k}(\lambda \xi)}{\gamma_{a,k}(\xi)} = 1
\,\,\,\,\textrm{for each}\,\, \lambda>0.$$ In other words,
see~\cite{karamata}, the function $\gamma_{a,k}$ is \textit{slowly
oscillating at infinity}. Moreover, in this case $\gamma_{a,k}$
may be represented by
$$\gamma_{a,k}(\xi) = \sigma_{a,k}(\xi)\,\exp\left(\int_1^\xi
\omega_{a,k}(t)\,\frac{dt}{t}\right),$$ where $\sigma_{a,k},
\omega_{a,k}\in C(\mathbb{R}_+)$ such that $$\lim_{\xi\to\infty}
\sigma_{a,k}(\xi)=\sigma>0, \qquad \lim_{\xi\to\infty}
\omega_{a,k}(\xi)=0.$$
\end{rem}

Similarly as in the case of algebras
$\mathcal{T}_k(\mathcal{A}_\infty)$ we get the following result
for algebras
$\mathcal{T}_k\left(L^{\{0,+\infty\}}_\infty(\mathbb{R}_+)\right)$.

\begin{thm}
Each $C^*$-algebra
$\mathcal{T}_k\left(L^{\{0,+\infty\}}_\infty(\mathbb{R}_+)\right)$
is isomorphic and isometric to $C[0,+\infty]$ and the isometric
isomorphism
$$\tau_k: \mathcal{T}_k\left(L^{\{0,+\infty\}}_\infty(\mathbb{R}_+)\right) \longrightarrow
C[0,+\infty]$$ is generated by the following mapping of generators
of
$\mathcal{T}_k\left(L^{\{0,+\infty\}}_\infty(\mathbb{R}_+)\right)$
$$\tau_k: T_a^{(k)} \longmapsto \gamma_{a,k}(\xi),$$
where $a=a(v)\in L^{\{0,+\infty\}}_\infty(\mathbb{R}_+)$.
\end{thm}

\noindent The inclusion
$$\tau_k\left(\mathcal{T}_k(L^{\{0,+\infty\}}_\infty(\mathbb{R}_+))\right)\subseteq C[0,+\infty]$$
is obvious from Theorem~\ref{thmlimits}, whereas the inverse
inclusion follows from the next theorem which is a consequence of
the Stone-Weierstrass theorem. Given a function $a_1(v)$ denote by
$L(1,a_1)$ the linear two-dimensional space generated by $1$ and
the function $a_1$.

\begin{thm}\label{thmgammaseparates}
Let $a_1(v)\in L^{\{0,+\infty\}}_\infty(\mathbb{R}_+)$ be a
real-valued function such that the corresponding function
$\gamma_{a_1,k}(\xi)$ separates the points of
$\overline{\mathbb{R}}_+$. Then each $C^*$-algebra
$\mathcal{T}_k(L(1,a_1))$ is isomorphic and isometric to
$C[0,+\infty]$ and the isometric isomorphism
$$\tau_k: \mathcal{T}_k(L(1,a_1)) \longrightarrow
C[0,+\infty]$$ is generated by the same mapping
$$\tau_k: T_a^{(k)} \longmapsto \gamma_{a,k}(\xi)$$
of generators of the algebra $\mathcal{T}_k(L(1,a_1))$.
\end{thm}

\begin{rem}\rm\label{rmkseparates}
As a consequence of this result we have that given a point
$\lambda_0\in\mathbb{R}_+$, then \textit{each $C^*$-algebra
$\mathcal{T}_k(L(1,\chi_{[0,\lambda_0]})$ is isomorphically
isometric to $C[0,+\infty]$}. In fact, by
Theorem~\ref{thmgammaseparates} we need to show that for each
$k\in\mathbb{Z}_+$ the real-valued function
$$\gamma_{\chi_{[0,\lambda_0]},k}(\xi) =
\chi_+(\xi)\int_{\mathbb{R}_+}
\chi_{[0,\lambda_0]}\left(\frac{v}{2\xi}\right)\ell_k^2(v)\,\mathrm{d}v
= \chi_+(\xi)\int_{0}^{2\lambda_0 \xi} \ell_k^2(v)\,\mathrm{d}v$$
separates the points of $\overline{\mathbb{R}}_+$. Instead of this
we have, see Appendix, that the function $\int_0^x
\Lambda_{0,k,k}^{(0,0)}(t)\,\mathrm{d}t$ is strictly increasing.
\end{rem}

The previous results motivate the study of piece-wise constant
symbols and algebras generated by Calder\'on-Toeplitz operators
with such symbols. But we will continue in the following
direction. Consider the Calder\'on-Toeplitz operator
$T_{a_+}^{(0)}$ with symbol $a_{+}(v)=\chi_{[0,1/2]}(v)$, which is
unitarily equivalent to $\gamma_{a_+,0} I$, where
$$\gamma_{a_+,0}(\xi) = 2\xi \int_{\mathbb{R}_+}
\chi_{[0,1/2]}(v)\, \mathrm{e}^{-2v\xi}\,\mathrm{d}v =
1-\mathrm{e}^{-\xi}, \quad \xi\in\overline{\mathbb{R}}_+.$$ This
function is continuous on $\overline{\mathbb{R}}_+$ with values in
$[0,1]$, therefore $T_{a_+}^{(0)}$ is self-adjoint and
$\textrm{sp}\, T_{a_+}^{(0)}=[0,1]$. Also, the function
$\gamma_{a_+,0}$ is strictly increasing and its inverse has the
form
$$\gamma_{a_+,0}^{-1}(x)=\xi(x)=-\ln(1-x),\quad x\in [0,1].$$ Thus,
for any function $h$ continuous on $[0,1]$ the operator
$h\left(T_{a_+}^{(0)}\right)$ is well defined according to the
standard functional calculus in $C^*$-algebras. Now we will
exploit the isomorphism between the Calder\'on-Toeplitz operator
algebra and the functional algebra given in
Corollary~\ref{coralgebra}.

\begin{thm}\label{thmexchange1}
Let $\Delta_\lambda$ be a family of functions parameterized by
$\lambda\in\mathbb{R}_+$ and given by
$$\Delta_\lambda(x)=1-(1-x)^{2\lambda}, \quad x\in [0,1].$$ Then $$\Delta_\lambda\left(T_{a_+}^{(0)}\right) =
T_{\chi_{[0,\lambda]}}^{(0)}\in
\mathcal{T}_0(L(1,\chi_{[0,\lambda]})).$$
\end{thm}

\begin{rem}\rm
Observe that each function $\Delta_\lambda$ is continuous on
$[0,1]$ and $\Delta_\lambda(0)=0$, $\Delta_\lambda(1)=1$ for each
$\lambda\in\mathbb{R}_+$. Some particular cases are also
interesting: $\Delta_{1/2}(x)\equiv x$ for all $x\in[0,1]$,
whereas for limit values as $\lambda\to 0$ and $\lambda\to+\infty$
we have
$$\Delta_0(x)\equiv 0 \quad \textrm{and}\quad
\Delta_\infty(x)\equiv 1,$$ respectively, for all $x\in (0,1)$.
These particular cases lead to the equalities
$$\Delta_{0}\left(T_{a_+}^{(0)}\right) = 0, \quad
\Delta_{1/2}\left(T_{a_+}^{(0)}\right) = T_{a_+}^{(0)}, \quad
\Delta_{\infty}\left(T_{a_+}^{(0)}\right) = I.$$
\end{rem}

We have chosen the Calder\'on-Toeplitz operator $T_{a_+}^{(0)}$ as
the starting operator because in this specific case the equation
$$x=\gamma_{a_+,0}(\xi) = 1-\mathrm{e}^{-\xi}$$ admits an explicit solution. But we can
start from any Calder\'on-Toeplitz operator
$T_{\chi_{[0,\lambda]}}^{(0)}$ with symbol
$a(v)=\chi_{[0,\lambda]}(v)$, $\lambda\in\mathbb{R}_+$. Indeed,
the function $\gamma_{\chi_{[0,\lambda]},0}(\xi)$ is strictly
increasing which implies that the function $\Delta_\lambda:
[0,1]\to [0,1]$ is strictly increasing as well and thus the
function $\Delta_\lambda^{-1}$ is well defined and continuous on
$[0,1]$. Clearly, for $\lambda_1, \lambda_2\in\mathbb{R}_+$ we
have $$\left(\Delta_{\lambda_2}\circ
\Delta_{\lambda_1}^{-1}\right)\left(T^{(0)}_{\chi_{[0,\lambda_1]}}\right)
=
T^{(0)}_{\chi_{[0,\lambda_2]}},$$ 
where $\circ$ is the usual composition of real functions. This
means that for any symbol $a=a(v)=\chi_{[0,\lambda]}(v)$ the
Calder\'on-Toeplitz operator $T_a^{(0)}$ belongs to the algebra
$\mathcal{T}_0(L(1,\chi_{[0,\lambda]}))$ and is the function of
the operator $T_{a_+}^{(0)}$, i.e.,
$$T_{\chi_{[0,\lambda]}}^{(0)}=\Delta_\lambda\left(T_{a_+}^{(0)}\right).$$ Moreover,
Theorem~\ref{thmgammaseparates} and Remark~\ref{rmkseparates}
imply that \textit{each Calder\'on-Toeplitz operator} with
$L_\infty^{\{0,+\infty\}}(\mathbb{R}_+)$-symbol can be obtained in
a similar way.

\begin{thm}\label{thmexchange2}
Let $a=a(v)\in L_\infty^{\{0,+\infty\}}(\mathbb{R}_+)$, and for
$x\in [0,1]$ put
$$\nabla_{a,\lambda}^{(k)}(x)=-\frac{1}{\lambda}\ln(1-x)\int_{\mathbb{R}_+} a(v)(1-x)^{v/\lambda}
L_k^2\left(-\frac{v}{\lambda}\ln(1-x)\right)\,\mathrm{d}v,$$ where
$\lambda\in\mathbb{R}_+$. Then
$$\nabla_{a,\lambda}^{(k)}\left(T_{\chi_{[0,\lambda]}}^{(0)}\right) =
T_{a}^{(k)}.$$
\end{thm}

\proof By Theorem~\ref{CTO1} the Calder\'on-Toeplitz operator
$T_{a}^{(k)}$ is unitarily equivalent to $\gamma_{a,k}I$, where
$$\gamma_{a,k}(\xi) = 2\xi
\int_{\mathbb{R}_+} a(v) \ell_k^2(2v\xi)\,\mathrm{d}v, \quad
\xi\in\overline{\mathbb{R}}_+.$$ Since
$x=\gamma_{\chi_{[0,\lambda]},0}(\xi)=1-\mathrm{e}^{-2\lambda\xi}$,
substituting $\xi=\xi(x)=-(2\lambda)^{-1}\ln(1-x)$ we have
$$\gamma_{a,k}(\xi(x)) = -\frac{1}{\lambda}\ln(1-x)
\int_{\mathbb{R}_+}
a(v)\ell_k^2\left(-\frac{v}{\lambda}\ln(1-x)\right)\,\mathrm{d}v =
\nabla_{a,\lambda}^{(k)}(x),$$ which completes the proof. \qed

\paragraph{Interpretations and applications} Consider an arbitrary
\textit{Toeplitz operator} acting on the Bergman space
$\mathcal{A}^2(\Pi)$ with symbol $\chi_{[0,\lambda]}(\Im \zeta)$
where $\zeta=u+\mathrm{i} v\in\Pi$. Using this operator
\textit{each Calder\'on-Toeplitz operator} with symbol $a=a(v)\in
L_\infty^{\{0,+\infty\}}(\mathbb{R}_+)$ acting on wavelet subspace
$A^{(k)}$ is the function of Toeplitz operator. 
This result is interesting because it enables to change not only
symbol from a "nice class" (as it is in the case of classical
Toeplitz operators), but also wavelet as will be explained now.

Let $k\in\mathbb{Z}_+$. In Section~\ref{section2} we have denoted
by $\mathcal{W}_k f$ the continuous wavelet transform of $f\in
H_2^+(\mathbb{R})$ with respect to wavelet $\psi^{(k)}$. Then the
value
$$(\mathcal{W}_k f)(\zeta) = \frac{1}{\sqrt{v}} \int_\mathbb{R}
f(x)\,\overline{\psi^{(k)}\left(\frac{x-u}{v}\right)}\,\mathrm{d}x,
\quad \zeta=(u,v)\in G,$$ is the wavelet coefficient at time $u$
and scale $v$. This integral measures the comparison of the local
shape of the signal $f$ and the shape of wavelet $\psi^{(k)}$.
Since $(\mathcal{W}_k f)(\zeta) \in A^{(k)}$, we say that this
comparison is made "on the level $k$". It is well-known that the
change of the value of the dilation factor $v$ serves as a
mathematical microscope to zoom in and out of the signal whereas
localization in time is achieved by selecting $u$. Thus, some time
and frequency localization is achieved for each point
$\zeta=(u,v)$ in the wavelet half-plane.

It is well-known that $(\mathcal{W}_k f)(\zeta)$ contains enough
information to reconstruct the function $f$ on the level $k$. This
inverse transform expresses the fact that no information is lost
in the transform and we have a representation of the signal $f$ on
the level $k$ as a linear superposition of wavelets
$\left(\rho_\zeta \psi^{(k)}\right)$ with coefficients
$(\mathcal{W}_k f)(\zeta)$. In fact, this "perfect reconstruction
of signal $f$ on level $k$" corresponds to the application of
Calder\'on-Toeplitz operator $T_a^{(k)}$ with symbol $a(u,v)\equiv
1$ to the signal $f$. Thus, emphasizing or eliminating some
information content in wavelet half-plane we get a filtered
version of signal $f$ for other choices of symbol $a$. Therefore,
these operators are a version of non-stationary (or, time-varying)
filters.

Restrict our attention to the case of symbols $a$ depending only
on vertical coordinate $v=\Im\zeta$, $\zeta\in G$, in wavelet
half-plane. Since $v$ is a measure of the duration of the event
being examined, then the operator $T^{(0)}_{a_+}$ gives a
reconstruction of a signal on the segment
$\Omega_{1/2}=\mathbb{R}\times (0,1/2]$ and level $0$. By the
result of Theorem~\ref{thmexchange1} we may obtain any operator
$T_{\chi_{[0,\lambda]}}^{(0)}$ (giving a reconstruction of a
signal on the segment $\Omega_\lambda=\mathbb{R}\times
(0,\lambda]$ and level $0$) from this operator $T_{a_+}^{(0)}$. In
fact, \textsf{from the reconstruction of a signal on the segment
$\Omega_{1/2}$ and level $0$ we may obtain a reconstruction of the
same signal on the same level on an arbitrary segment
$\Omega_\lambda$} using the function $\Delta_\lambda$ which is
easy to compute.

Even more interesting is the result of Theorem~\ref{thmexchange2}.
If we know the reconstruction of a signal on a segment
$\Omega_\lambda$ and level $0$, we might get an arbitrary
reconstruction of the signal (as its filtered version using a real
bounded function $a$ of scale having limits in critical points of
boundary of $\mathbb{R}_+$ such that the corresponding function
$\gamma_\cdot$ separates the points of $\overline{\mathbb{R}}_+$)
on an arbitrary level $k$ using the function
$\nabla_{a,\lambda}^{(k)}$. Theoretically, for the purpose to
study localization of a signal in the wavelet half-plane the
result of Theorem~\ref{thmexchange2} suggests to consider certain
"nice" symbols on the first level $0$ (indeed, Toeplitz operators
on $\mathcal{A}^2(\Pi)$ with symbols as characteristic functions
of some interval in $\mathbb{R}_+$) instead of possibly
complicated $L_\infty^{\{0,+\infty\}}(\mathbb{R}_+)$-symbols with
respect to "different microscope" represented by the level $k$. On
the other hand, to compute the corresponding function
$\nabla_{a,\lambda}^{(k)}$ need not be always easy.

\section{Appendix: functions and integrals containing Laguerre polynomials}\label{sectionE}

Here we state some necessary formulas and estimations for
functions and integrals containing Laguerre polynomials which were
used on different places of this article. First we recall the
following important integral formula, cf.~\cite[formula (16),
p.~330]{WG},
\begin{align}\label{Laguerreformula}
& \phantom{=} \int_{\mathbb{R}_+} x^{p} \mathrm{e}^{-x}
L_m^{(\alpha)}(x)L_n^{(\beta)}(x)\,\mathrm{d}x \\ & = \Gamma(p+1)
\sum_{i=0}^{\min\{m,n\}} (-1)^{m+n}{p-\alpha \choose m-i} {p-\beta
\choose n-i} {p+i \choose i},\nonumber
\end{align} where $\Re p>-1$, $\alpha,\beta>-1$,
$m,n\in\mathbb{Z}_+$, and $${a \choose b} =
\frac{\Gamma(a+1)}{\Gamma(b+1) \Gamma(a-b+1)}.$$ Applying the
formula~(\ref{Laguerreformula}) for $p=1$, $\alpha=\beta=0$ and
$m=n=k$ we get the result in formula~(\ref{lambda_k}). Also, the
formula~(\ref{Laguerreformula}) plays an important role in many
areas of research, e.g., in mathematical physics (in the quantum
mechanical treatment of the hydrogen atom).

Recall that the derivative of generalized Laguerre polynomial is
given by the formula
\begin{equation}\label{derivativeL_k}
\frac{d^r}{dx^r} L_k^{(\alpha)}(x) = (-1)^r
L_{k-r}^{(\alpha+r)}(x),
\end{equation} cf.~\cite[formula 8.971.2]{GR}. Note that for $k<r$ the $r$-th derivative of $L_k^{(\alpha)}$ is always
zero, and thus we always use the convention $L^{(\alpha)}_n(x)=0$
for $n<0$ and arbitrary $\alpha>-1$. The
formula~(\ref{derivativeL_k}) together with the Leibniz rule
yields the $n$-th derivative of function $\ell_k^2(2v\xi)$, i.e.,
\begin{align}\label{n-thderivativeell}
\frac{\mathrm{d}^n}{\mathrm{d}\xi^n}\ell_k^2(2v\xi) & =
\sum_{i=0}^n {n\choose
i}\frac{\mathrm{d}^{n-i}}{\mathrm{d}\xi^{n-i}} \mathrm{e}^{-2v\xi}
\cdot\sum_{j=0}^i {i\choose j}
\frac{\mathrm{d}^{i-j}}{\mathrm{d}\xi^{i-j}}
L_k(2v\xi) \frac{\mathrm{d}^{j}}{\mathrm{d}\xi^{j}} L_k(2v\xi) \nonumber \\
& = (-2v)^{n} \mathrm{e}^{-2v\xi}\sum_{i=0}^n \sum_{j=0}^i {n
\choose i} {i \choose j} L_{k-i+j}^{(i-j)}(2v\xi)
L_{k-j}^{(j)}(2v\xi).
\end{align}

On many places in this paper the estimation of the following
(non-negative) function $$\Lambda_{p,m,n}^{(\alpha,\beta)}(x) :=
x^p \mathrm{e}^{-x} \left|L_m^{(\alpha)}(x)
L_n^{(\beta)}(x)\right|, \quad x\in\mathbb{R}_+,$$ has been used.
In our case all parameters $p, m, n, \alpha,
\beta\in\mathbb{Z}_+$, but some results hold also for $p, \alpha,
\beta$ suitable real numbers, even $p$ may be a suitable complex
number with $\Re p>-1$. Clearly, for all admissible values of
parameters the function
\begin{equation}\label{intLambda}
\int_0^x \Lambda_{p,m,n}^{(\alpha,\beta)}(t)\,\mathrm{d}t, \quad
x\in\mathbb{R}_+,\end{equation} is non-negative, continuous,
strictly increasing and bounded on $\mathbb{R}_+$. For the latter
property see~(\ref{upperbound2}). For the purpose of boundedness
of~(\ref{intLambda}) we need the following result. Recall that
$$(x)_n = \frac{\Gamma(x+n)}{\Gamma(x)}$$ is the Pochhammer
symbol.

\begin{thm}[\cite{lewandowski}]
For $\alpha\geq -\frac{1}{2}$, $x\geq 0$ and $n\in\mathbb{Z}_+$ we
have
$$\left|L_n^{(\alpha)}(x)\right| \leq \frac{(\alpha+1)_n}{n!}\,
\sigma_n^{(\alpha)}(\mathrm{e}^x),$$ where $\sigma_n^{(\alpha)}$
are the Ces\`{a}ro means of the formal series
$\sum\limits_{i=0}^\infty b_i$ given by the formula
$$\sigma_n^{(\alpha)}\left(\sum_{i=0}^\infty b_i\right) := \frac{n!}{(\alpha+1)_n}
\sum_{i=0}^n \frac{(\alpha+1)_{n-i}}{(n-i)!}\,b_i, \quad
\alpha>-1.$$
\end{thm}
Using this theorem, for each $\alpha,\beta\geq -\frac{1}{2}$,
$p>-1$, $x\in\mathbb{R}_+$ and $m,n\in\mathbb{Z}_+$ we immediately
have
\begin{align}\label{upperbound1}
\Lambda_{p,m,n}^{(\alpha,\beta)}(x) & \leq \frac{(\alpha+1)_m}{m!}
\frac{(\beta+1)_n}{n!} x^p\, \mathrm{e}^{-x}
\sigma_{m}^{(\alpha)}(\mathrm{e}^x) \sigma_{n}^{(\beta)}(\mathrm{e}^x) \nonumber \\
& = x^p\, \mathrm{e}^{-x} \sum_{i=0}^m
\frac{(\alpha+1)_{m-i}}{(m-i)!\,i!}\,x^i\cdot \sum_{j=0}^n
\frac{(\beta+1)_{n-j}}{(n-j)!\,j!}\,x^j \nonumber \\ & =
\sum_{i=0}^m \sum_{j=0}^n \frac{(\alpha+1)_{m-i}}{(m-i)!\,i!}
\frac{(\beta+1)_{n-j}}{(n-j)!\,j!}\,x^{p+i+j}\, \mathrm{e}^{-x},
\end{align} and thus for its integral holds
\begin{align}\label{upperbound2}
\int_{\mathbb{R}_+}
\Lambda_{p,m,n}^{(\alpha,\beta)}(x)\,\mathrm{d}x & \leq
\sum_{i=0}^m \sum_{j=0}^n \frac{(\alpha+1)_{m-i}}{(m-i)!\,i!}
\frac{(\beta+1)_{n-j}}{(n-j)!\,j!}\,\int_{\mathbb{R}_+}
x^{p+i+j}\, \mathrm{e}^{-x}\,\mathrm{d}x \nonumber \\ & =
\sum_{i=0}^m \sum_{j=0}^n \frac{(\alpha+1)_{m-i}}{(m-i)!\,i!}
\frac{(\beta+1)_{n-j}}{(n-j)!\,j!}\,\Gamma(p+i+j+1) \nonumber \\ &
:= \textrm{const}_{p,m,n}^{(\alpha,\beta)}.
\end{align}

The following observation is easy to verify using the classical
methods of calculus: for each $p>0$, $q>0$ and each
$x\in\mathbb{R}_+$ holds
\begin{equation}\label{estimation}
x^p \mathrm{e}^{-q x} \leq \left(\frac{p}{\mathrm{e}\,
q}\right)^p.
\end{equation}

The following result was very useful in proof of
Theorem~\ref{thmlimits}. It may be of further interest itself,
because in some sense it generalizes the incomplete Gamma function
to general Laguerre functions. As far as we know it is not
included in any literature we have seen, so for the sake of
completeness we give its short proof here which is only of
computational nature.

\begin{lem}\label{lemmaell}
For each $k\in\mathbb{Z}_+$ and each $x\in\mathbb{R}_+$ we have
$$\int_0^x \Lambda_{0,k,k}^{(0,0)}(t)\,\mathrm{d}t = \int_0^x \ell_k^2(t)\,\mathrm{d}t
= 1-N_{2k}(x)\,\mathrm{e}^{-x},$$ where $N_{2k}(x)$ is a
polynomial of $x$ of degree $2k$.
\end{lem}

\proof Integrating by parts we have $$\int_0^x
\Lambda_{0,k,k}^{(0,0)}(t)\,\mathrm{d}t = \int_0^x \mathrm{e}^{-t}
L_k^2(t)\,\mathrm{d}t = 1-\mathrm{e}^{-x} L_k^2(x) -2\int_0^x
\mathrm{e}^{-t}L_k(t) L_{k-1}^{(1)}(t)\,\mathrm{d}t,
$$ where the formula~(\ref{derivativeL_k}) has been used. The
last integral has its sense for $k\geq 1$, whereas for $k=0$
disappears and we have immediately the desired result. Using the
explicit form of Laguerre polynomials $L^{(\alpha)}_n$ given
in~(\ref{defLaguerre}) we get
$$\int_0^x \mathrm{e}^{-t}L_k(t) L_{k-1}^{(1)}(t)\,\mathrm{d}t =
\sum_{i=0}^k \sum_{j=0}^{k-1} (-1)^{i+j} {k\choose i}{k\choose
j+1} \frac{1}{i! j!} \int_0^x \mathrm{e}^{-t}
t^{i+j}\,\mathrm{d}t.$$ The last integral is, in fact, the
incomplete Gamma function $\gamma(1+i+j,x)$ for which the formula
$$\gamma(1+n,x) = n!\left(1-\mathrm{e}^{-x}
\sum_{p=0}^{n} \frac{x^{p}}{p!}\right), \quad n\in\mathbb{Z}_+,$$
holds, see~\cite[formula 8.352.1]{GR}. Summarizing the above
yields
\begin{align*}
& \phantom{=} \int_0^x \Lambda_{0,k,k}^{(0,0)}(t)\,\mathrm{d}t =
1-2\sum_{i=0}^k
\sum_{j=0}^{k-1}(-1)^{i+j} {k\choose i}{k\choose j+1} {i+j\choose i} \\
& - \mathrm{e}^{-x}\left(L_k^2(x)-2\sum_{i=0}^k
\sum_{j=0}^{k-1}\sum_{p=0}^{i+j} (-1)^{i+j} {k\choose i}{k\choose
j+1}{i+j\choose i} \frac{x^{p}}{p!}\right).
\end{align*}Since the expression in brackets is a polynomial of $x$ of
degree $2k$, it suffices to show that $$S(k) := \sum_{i=0}^k
\sum_{j=0}^{k-1}(-1)^{i+j} {k\choose i}{k\choose j+1} {i+j\choose
i} = 0$$ for each $k\geq 1$. But for each $k\geq 1$ we have $$S(k)
= \sum_{i=0}^k (-1)^{i} {k\choose i} \cdot \sum_{j=0}^{k-1}
\frac{(-1)^{j}}{j!}{k\choose j+1} (i+1)_j = \sum_{i=0}^k (-1)^{i}
{k\choose i} Q(i),$$ where $Q(x)$ is a polynomial of degree $k-1$.
Then the result from the theory of finite differences yields
$S(k)=0$. \qed

\vspace{12pt} {\bf Acknowledgements:} Author has been on a
postdoctoral stay at the Departamento de Matem\'a\-ti\-cas,
CINVESTAV del IPN (M\'exico), when writing this paper and
investigating the topics presented herein. He therefore gratefully
acknowledges the hospitality and support of the mathematics
department of CINVESTAV on this occasion. Especially, author
wishes to thank Professor Nikolai L. Vasilevski for giving him the
book~\cite{vasilevskibook} and for his postdoctoral guidance and
help during his stay in M\'exico. Also, autor's thanks go to
Professor Nico M. Temme for drawing his attention to the
paper~\cite{lewandowski} with useful boundedness results for
Laguerre polynomials. Moreover, author would like to express his
deepest gratitude to anonymous referees for their suggestions and
interesting observations. Especially, many important consequences
of Theorem~\ref{thmlimits} presented in Introduction are due to
the anonymous referee, and we appreciate them very well.

\vspace{5mm}

\noindent \small{Ondrej Hutn\'ik, Departamento de Matem\'aticas,
CINVESTAV del IPN, {\it Current address:} Apartado Postal 14-740,
07000, M\'exico, D.F., M\'exico
\newline {\it E-mail address:} hutnik@math.cinvestav.mx}
\newline AND \newline \noindent \small{Institute of Mathematics,
Faculty of Science, Pavol Jozef \v Saf\'arik University in Ko\v
sice, Jesenn\'a 5, 040~01 Ko\v sice, Slovakia,
\newline {\it E-mail address:} ondrej.hutnik@upjs.sk}

\end{document}